\documentclass[reqno,twoside,12pt,a4paper]{amsart}
\usepackage{latexsym}
\usepackage{amsmath}
\usepackage{amssymb}
\usepackage{cases}
\usepackage{mathrsfs}
\usepackage{amsthm}
\usepackage{cite}
\usepackage[usenames,dvipsnames]{color}

\def\pier #1{{\color{red}#1}}
\def\takeshi #1{{\color{blue}#1}}
\def\hao #1{{\color{green}#1}}
\def\pcol #1{{\color{red}#1}}
\def\pier #1{#1}
\def\takeshi #1{#1}
\def\hao #1{#1}
\def\pcol #1{#1}

\allowdisplaybreaks[4]

\title[Transmission problem of {C}ahn--{H}illiard type]{On a transmission problem for equation and \\ dynamic boundary condition of {C}ahn--{H}illiard type\\ with nonsmooth potentials}

\author[P. Colli]{Pierluigi Colli}
\address{Pierluigi Colli: 
Dipartimento di Matematica, Universit\`a degli Studi di Pavia\\
\pcol{and Research Associate at the IMATI -- C.N.R. Pavia}\\
Via Ferrata~\pcol{5}, 27100 Pavia, Italy}
\email{pierluigi.colli@unipv.it}

\author[T.\ Fukao]{Takeshi Fukao}
\address{Takeshi Fukao: Department of Mathematics, Faculty of Education, 
Kyoto University of Education, 
1~Fujinomori, Fukakusa, Fushimi-ku, Kyoto~612-8522 Japan}
\email{fukao@kyokyo-u.ac.jp}

\author[H.\ Wu]{Hao Wu}
\address{Hao Wu: School of Mathematical Sciences and Shanghai Key Laboratory for Contemporary Applied Mathematics, Fudan University, Han Dan Road 220, Shanghai 200433, China; Key Laboratory of Mathematics for Nonlinear Science (Fudan University),
Ministry of Education, Han Dan Road 220, Shanghai 200433, China}
\email{haowufd@fudan.edu.cn}

\dedicatory{}
\subjclass[2000]{}
\begin{document}

\thispagestyle{empty}
\maketitle

\begin{abstract}
This paper is concerned with well-posedness \pier{of} 
the {C}ahn--{H}illiard equation subject to a class of new dynamic boundary conditions. 
The system was recently derived in \takeshi{Liu--Wu \pcol{(Arch. Ration. Mech. Anal. {\bf 233} (2019), 167--247)}} via an energetic variational approach and it naturally fulfills three physical constraints such as mass conservation, energy dissipation and force balance. The target problem \pier{examined in} this paper can be viewed as a transmission problem that consists \pier{of} {C}ahn--{H}illiard type equations both in the bulk and on the boundary\pier{.} \pier{In our approach, we are able to deal with a general class of potentials with double-well structure, including the physically relevant logarithmic potential and the non-smooth double-obstacle potential. Existence, uniqueness and continuous dependence of global weak solutions are established. The proof is based on a novel 
time-discretization scheme for the approximation of the continuous problem. Besides, a regularity result is shown with the aim of obtaining a strong solution to the system.} 
\vskip3mm
\pier{\noindent {\sc Key words:}
Cahn--Hilliard system, dynamic boundary condition, transmission problem, 
non-smooth potentials, well-posedness, regularity.
}
\end{abstract}

%%%%% Section 1. %%%%%
\section{Introduction}
\setcounter{equation}{0}

In this paper, we consider the following initial boundary value problem for a {C}ahn--{H}illiard equation subject to a dynamic boundary condition that is also of {C}ahn--{H}illiard type. Let $0<T<\infty$ be some fixed time and let $\Omega \subset \mathbb{R}^{d}$, 
$d=2$ or $3$, be a bounded domain with smooth boundary $\Gamma := \partial \Omega$. \pier{We aim to find} four unknown functions 
$\phi, \mu: Q:=(0,T) \times \Omega \to \mathbb{R}$ and 
$\psi, w: \Sigma:=(0,T) \times \Gamma$ \pier{satisfying}
\begin{align} 
	&\partial_t \phi - \Delta \mu =0 
	& \mbox{in }Q, 
	\label{lw1}
	\\
	&\mu = -\Delta \phi + {\mathcal W}'(\phi) 
	& \mbox{in }Q, 
	\label{lw2}
	\\
	&\partial_{\boldsymbol{\nu}} \mu =0 
	& \mbox{on }\Sigma, 
	\label{lw3}
	\\
	&\phi_{|_{\Gamma}} =\psi 
	& \mbox{on }\Sigma,
	\label{lw4}
	\\ 
	&\partial_t \psi -\Delta_\Gamma w = 0 
	& \mbox{on }\Sigma, 
	\label{lw5}
	\\
	&w = \partial _{\boldsymbol{\nu}} \phi - \Delta_\Gamma \psi 
	+ {\mathcal W}_\Gamma'(\psi)
	& \mbox{on }\Sigma,
	\label{lw6}
\end{align}
where $\partial_t$ and $\partial_{\boldsymbol{\nu}}$ denote 
the partial time derivative and the outward normal derivative \pier{on $\Gamma,$} respectively; 
$\Delta$ denotes the {L}aplacian and  
$\Delta _{\Gamma }$ denotes the {L}aplace--{B}eltrami operator
on $\Gamma $ (see, e.g., \cite[Chapter~3]{Gri09}); 
$\phi_{|_\Gamma}$ standards for the trace of $\phi $ on the boundary $\Gamma$. In view of \eqref{lw4}, system (\ref{lw1})--(\ref{lw6}) is a sort of transmission problem between the {C}ahn--{H}illiard equation in the bulk 
$\Omega$ and the {C}ahn--{H}illiard equation on the boundary $\Gamma$.  
\pier{The nonlinear functions ${\mathcal W}$ and ${\mathcal W}_\Gamma$ are usually referred 
as the double-well potentials, with two minima and a local unstable maximum in between. 
Typical and physically significant examples of such potentials 
are the so-called {\em classical potential}, the {\em logarithmic potential\/},
and the {\em double obstacle potential\/}, which are given, in this order,~by
\begin{align}
  & {\mathcal W}_{reg}(r) := \frac 14 \, (r^2-1)^2 \,,
  \quad r \in \mathbb{R}, 
  \label{regpot}
  \\
  & {\mathcal W}_{log}(r) := (1+r)\ln (1+r)+(1-r)\ln (1-r)  - c_1 r^2 \,,
  \quad r \in (-1,1),
  \label{logpot}
  \\[1mm]
  & {\mathcal W}_{2obs}(r) := c_2 (1-r^2) 
  \quad \hbox{if $|r|\leq1$}
  \quad\hbox{and} \quad
  {\mathcal W}_{2obs}(r) := +\infty
  \quad \hbox{if $|r|>1$}.
  \label{obspot}
\end{align}
where the constants in \eqref{logpot} and \eqref{obspot} satisfy
$c_1>1$ and $c_2>0$, so that
${\mathcal W}_{log}$ and ${\mathcal W}_{2obs}$ are nonconvex.}
\pier{The nonlinear terms ${\mathcal W}'(\phi) $ in \eqref{lw2} and ${\mathcal W}'_\Gamma (\psi)
$ in \eqref{lw6} \pier{characterize} the dynamics of the {C}ahn--{H}illiard system. In cases like 
\eqref{regpot} and \eqref{logpot}, ${\mathcal W}'$ and ${\mathcal W}'_\Gamma$
denote simply the derivatives of the related potentials; while non-smooth potentials like \eqref{obspot} are considered, then ${\mathcal W}'$ and ${\mathcal W}'_\Gamma$
denote the subdifferential of the convex part plus the derivative of the smooth concave contribution, i.e., for \eqref{obspot} it is 
$$
s\in {\mathcal W}'_{2obs}(r) \quad \hbox{if} \quad r\in [-1,1], \quad s +2  c_2 r \, 
\begin{cases}
\in (- \infty , 0 ] & \hbox{if } \ r=-1 \\ 
=0  & \hbox{if } \ r\in (-1, 1) \\ 
\in [0, + \infty ) &\hbox{if } \  r=1
\end{cases} . 
$$
Of course, in this case one should replace the equalities in \eqref{lw2} and \eqref{lw6}
by inclusions. In this paper, we are able to handle completely general potentials 
${\mathcal W}$ and ${\mathcal W}_\Gamma$ including all the three cases \eqref{regpot}--\eqref{obspot} mentioned above.}

The system (\ref{lw1})--(\ref{lw6}) was first derived by {L}iu and {W}u \cite{LW19} in a more general form (see also \cite{Wu18}) on the basis of an energetic variational approach. It describes effective short-range interactions between the binary mixture and the solid wall (boundary), furthermore, it has the feature that the related model naturally fulfills important physical constraints such as conservation of mass, dissipation of energy and force balance relations. In its current formulation, we see that 
\pier{equations (\ref{lw1}) and (\ref{lw2}) yield} a {C}ahn--{H}illiard system subject to a no-flux boundary condition (\ref{lw3}) together with a non-homogeneous Dirichlet boundary condition (\ref{lw4}), while the dynamic boundary condition (\ref{lw5}) and equation (\ref{lw6}) provide an evolution system of {C}ahn--{H}illiard type  on the boundary $\Gamma$.
These two {C}ahn--{H}illiard systems in the bulk and on the boundary 
are coupled through the trace condition (\ref{lw4}) and the normal derivative term $\partial_{\boldsymbol{\nu}} \phi$ in (\ref{lw6}).

\pier{The total energy functional for system (\ref{lw1})--(\ref{lw6}) given by 
\begin{equation}
	E(\phi, \psi):=\int_{\Omega} \takeshi{\left(  \frac{1}{2}|\nabla \phi|^2 + {\mathcal W}(\phi) \right) } dx
	+\int_{\Gamma} \takeshi{ \left( \frac{1}{2}|\nabla_\Gamma \psi|^2 + {\mathcal W}_\Gamma (\psi) \right) }d\Gamma
	\label{ene}
\end{equation}
is decreasing in time (see \cite{LW19}) and furthermore, system (\ref{lw1})--(\ref{lw6}) can be interpreted as a gradient flow of  $E(\phi, \psi)$ in a suitable dual space (see \cite{GK18}). In light of (\ref{lw1}), (\ref{lw3}) 
and (\ref{lw5}), we easily deduce that the following properties on mass conservation:
\begin{equation}
	\int_{\Omega} \phi(t) dx = \int_{\Omega} \phi(0) dx, 
	\quad  
	\int_{\Gamma} \psi(t) d\takeshi{\Gamma} = \int_{\Gamma} \psi(0) d\takeshi{\Gamma}
	\quad \mbox{for all } \, t \in [0,T]. 
	\label{mass}
\end{equation}
}

In this paper, we study the 
well-posedness of system (\ref{lw1})--(\ref{lw6}) for a weak solution subject to the following initial data
\begin{equation} 
	\phi(0)=\phi_0
	\quad 
	\mbox{in }\Omega, 
	\quad 
	\psi(0)=\psi_0
	\quad \mbox{on }\Gamma.
	\label{lw7}
\end{equation}
Moreover, we also \pier{establish a regularity theory in order to obtain} a strong solution. In particular, we are able to treat \pier{the initial value problem for system (\ref{lw1})--(\ref{lw6})
in a wider class of nonlinearities ${\mathcal W}$ and ${\mathcal W}_\Gamma$. 
Indeed, in the previous contributions, the well-posedness was investigated 
only in the case of smooth potentials like \eqref{regpot} (cf.~\cite[Remark~3.2]{LW19} and \cite[Remark 2.1]{GK18}): this is the point of emphasis of our present paper.}  

\pier{We would like to mention some related problems in the literature.
In 2011, {G}oldstein, {M}iranville and {S}chimperna~\cite{GMS11} studied a 
different type of transmission problem between the {C}ahn--{H}illiard system in the bulk and on the boundary with non-permeable walls (cf. a previous work   Gal \cite{Gal06} for the case with permeable walls). Their system can be derived 
from the same energy functional (\ref{ene}) by a variational method, however, the corresponding boundary conditions turn out to be different from (\ref{lw3}) and (\ref{lw5}). This also leads to a different property on the mass conservation comparing with (\ref{mass}) such that the total (bulk plus boundary) mass is conserved. 
We refer to \cite{LW19} for more detailed information on the comparison between these models. In addition, we mention the contributions \cite{CGM13, CF15, Gal06, GMS11} related to the well-posedness, \cite{FY17, FW19, CGS18, GS19,Gal08,GW08} for the study of long time behavior and the optimal control problems, \cite{CP14, FYW17} for numerical analysis and \cite{Kaj18} for 
the maximal regularity theory. 
Comparing the large number of known results on the previous model \cite{GMS11,Gal06}, 
we are only aware of the recent papers \cite{GK18, LW19} that analyze the well-posedness of system (\ref{lw1})--(\ref{lw6}) with (\ref{lw7}).}

\pier{Let us now describe the contents of the present paper.
In Section~2, we state the main well-posedness result for global weak solutions. 
We consider the problem within a general framework by setting  
${\mathcal W}':=\beta+\pi$ and ${\mathcal W}'_\Gamma:=\beta_\Gamma+\pi_\Gamma$, 
where $\beta$ and $\beta_\Gamma$ are maximal monotone graphs 
with $0\in \beta (0) $ and $0\in \beta_\Gamma (0) $, while 
$\pi$ and $\pi_\Gamma$ yield the anti-monotone terms that  
are {L}ipschitz continuous functions. 
The main theorems are concerned with the existence of a global weak solution (Theorem~2.1)  and the continuous dependence on the given data (Theorem~2.2), which implies the uniqueness.}

\pier{In Section~3, we study the time-discrete 
approximate problem for (\ref{lw1})--(\ref{lw6}) with (\ref{lw7}). 
We start from the viscous {C}ahn--{H}illiard system 
by inserting two additional terms, $\tau \partial_t \phi$ and $\sigma \partial_t \psi$
in the right hand sides of (\ref{lw2}) and (\ref{lw6}), respectively, with the 
parameters $\tau, \sigma>0$. Moreover, we take the {Y}osida approximations $\beta_\varepsilon$ and $\beta_{\Gamma, \varepsilon} $
in place of the maximal monotone graphs $\beta$ and $\beta_\Gamma$ and in terms 
of the parameter $\varepsilon>0$.  
Then we apply a time discretization scheme using the approach in \cite{CGS20, CK19}. We can show the existence of a discrete solution taking advantage of the general maximal monotone theory. After that, we proceed to derive a sequence of uniform estimates. For this purpose, we apply the technique of \cite{CC13} in order to treat different potentials in the bulk and on the boundary. In the subsequent iterations, we prove the existence results by performing the limiting procedures, 
with respect to the time step first, then as $\varepsilon \to 0$, finally 
taking the limit as either $\tau \to 0$ or $\sigma \to 0$, or both $\tau, \sigma \to 0$, in  order to obtain a partially viscous {C}ahn--{H}illiard system or a pure {C}ahn--{H}illiard system in the limit.
The continuous dependence result is then proved by using the energy method.}

\pier{In Section~4, we discuss the regularity for weak solutions. 
Returning to the time discrete approximation, we gain some necessary higher order estimates at all the different levels up to the final limits. 
Thus, we are able to obtain enough regularity as to guarantee a strong solution for the pure {C}ahn--{H}illiard system as well (see Theorem~4.1).}

\pier{Here, for the reader's convenience, let us include a detailed index of sections and subsections.}
\begin{itemize}
 \item[1.] Introduction
 \item[2.] Main results
 \item[3.] Well-posedness
 \begin{itemize}
  \item[3.1.] Time-discrete approximate solution
  \item[3.2.] A priori estimates and limiting procedure
  \item[3.3.] From viscous to pure {C}ahn--{H}illiard system
 \end{itemize}
 \item[4.] Existence of strong solution
\end{itemize}

%%%%% Section 2. %%%%%
\section{Main results}
\setcounter{equation}{0}

We now formulate our target problem \eqref{lw1}--\eqref{lw6} and \eqref{lw7} as follows: 
\begin{align} 
	&\partial_t \phi - \Delta \mu =0 
	& \mbox{a.e.\ in }Q, 
	\label{LW1}
	\\
	&\mu = -\Delta \phi + \xi + \pi (\phi) -f,\quad \xi \in \beta (\phi)
	& \mbox{a.e.\ in }Q, 
	\label{LW2}
	\\
	&\partial_{\boldsymbol{\nu}} \mu =0 
	& \mbox{a.e.\ on }\Sigma, 
	\label{LW3}
	\\
	&\phi_{|_{\Gamma}} =\psi 
	& \mbox{a.e.\ on }\Sigma,
	\label{LW4}
	\\ 
	&\partial_t \psi -\Delta_\Gamma w = 0 
	& \mbox{a.e.\ on }\Sigma, 
	\label{LW5}
	\\
	&	w = \partial _{\boldsymbol{\nu}} \phi - \Delta_\Gamma \psi 
	+ \zeta + \pi_\Gamma(\psi)-g,\quad 
	\zeta \in \beta_\Gamma (\psi)
	& \mbox{a.e.\ on }\Sigma,
	\label{LW6}
	\\
	&\phi(0)=\phi_0
	\quad 
	\mbox{a.e.\ in }\Omega, 
	\quad 
	\psi (0)=\psi_0
	\quad 
	 \mbox{a.e.\ on }\Gamma.
	 &
	\label{LW7}
\end{align}
where 
$f : Q \to \mathbb{R}$, $g : \Sigma \to \mathbb{R}$, 
$\phi_0:\Omega \to \mathbb{R}$, $\psi_0:\pier{\Gamma} \to \mathbb{R}$
 are given functions. 
Moreover, $\beta$ stands for the subdifferential of the convex part $\widehat{\beta }$ and 
$\pi $ stands for the derivative of the concave perturbation $\widehat{\pi}$ of a 
double well potential ${\mathcal W}(r)=\widehat{\beta }(r)+\widehat{\pi}(r)$. 
The same setting holds for $\beta _\Gamma $ and $\pi _\Gamma $.  
\pier{Typical examples of $\beta$, $\pi $ are given by (cf.~\eqref{regpot}--\eqref{obspot}):
\begin{itemize}
\item $\beta(r)=r^3$, $\pi (r)=-r$, $r\in \mathbb{R}$,  
for the prototype potential ${\mathcal W}_{reg}(r)$;\\
\item $\beta(r)=\ln((1+r)/(1-r))$, $\pi (r)=-2c_1r$, with 
$r \in (-1,1)$ for the logarithmic potential ${\mathcal W}_{log}(r)$;\\
\item $\beta(r)=\partial I_{[-1,1]}(r)$, $\pi (r)=-2c_2 r$, with 
$r \in [-1,1]$,  for the nonsmooth potential \takeshi{${\mathcal W}_{2obs}(r)$}.
\end{itemize}
Same considerations apply to $\beta_\Gamma$, $ \pi_\Gamma$ and ${\mathcal W}_\Gamma$. Since the bulk and boundary potentials are allowed to be different, in order to handle the nontrivial bulk-boundary interaction of the transmission problem, an assumption for the relationship between $\beta$ and $\beta_\Gamma$ will be needed. We shall present it later.} \smallskip

Hereafter, we use the spaces 
$$H:=L^2(\Omega ), \quad 
H_\Gamma :=L^2(\Gamma ), \quad 
V:=H^1(\Omega ), \quad 
V_\Gamma :=H^1(\Gamma )$$ 
with their dual spaces 
$V^*$ and $V_\Gamma^*$ of $V$ and $V_\Gamma$, respectively; and 
$$
W:=\{ z \in H^2(\Omega) : \partial _{\boldsymbol{\nu}} z =0 \ \mbox{a.e.\ on }\Gamma \}
$$   
equipped with the usual norms and inner products, 
denote them by  
$| \cdot |_{H}$ and $(\cdot,\cdot )_{H}$, and so on.\smallskip

Now, we define the weak solution of problem \eqref{LW1}--\eqref{LW7}:
\smallskip

\paragraph{\bf Definition~2.1.} 
{\it The sextuplet $(\phi, \mu, \xi, \psi, w, \zeta)$ is called \pier{a} weak solution of problem \eqref{LW1}--\eqref{LW7}, if
\begin{align*}
	& \phi \in H^1(0,T;V^*)\cap L^\infty (0,T;V) \cap L^2\bigl( 0,T;H^2(\Omega) \bigr), \\
	& \mu \in L^2(0,T;V), \quad \xi \in L^2(0,T;H), \\
	& \psi \in H^1(0,T;V_\Gamma^*)\cap L^\infty (0,T;V_\Gamma) \cap L^2 \bigl( 0,T;H^2(\Gamma) \bigr), \\
	& w \in L^2(0,T;V_\Gamma), \quad \zeta \in L^2(0,T;H_\Gamma)
\end{align*} 
and \pier{$\phi, \mu, \xi, \psi, w, \zeta$} satisfy
\begin{align} 
	& \langle \partial_t \phi, z \rangle_{V^*,V}+ \int_\Omega \nabla \mu \cdot \nabla z dx=0 \quad 
	\mbox{for all } z \in V, 
	& \mbox{a.e.\ in } (0,T),
	\label{weak1}
	\\
	& \mu = -\Delta \phi + \xi + \pi (\phi) -f,\quad \xi \in \beta (\phi)
	& \mbox{a.e.\ in }Q, 
	\label{weak2}
	\\
	& \phi_{|_{\Gamma}} =\psi 
	& \mbox{a.e.\ on }\Sigma,
	\label{weak3}
	\\ 
	& \langle \partial_t \psi,z_\Gamma  \rangle_{V_\Gamma^*,V_\Gamma}
	+ \int_{\Gamma} \nabla_\Gamma w \cdot \nabla_\Gamma z_\Gamma d\Gamma = 0 \quad 
	\mbox{for all } z_\Gamma \in V_\Gamma, 
	& \mbox{a.e.\ in } (0,T),
	\label{weak4}
	\\
	& 	w = \partial _{\boldsymbol{\nu}} \phi - \Delta_\Gamma \psi 
	+ \zeta + \pi_\Gamma(\psi)-g,  \quad \zeta \in \beta_\Gamma (\psi),
	& \mbox{a.e.\ on }\Sigma, 
	\label{weak5}
	\\
	&\phi(0)=\phi_0
	\quad 
	\mbox{a.e.\ in }\Omega, 
	\quad 
	\psi (0)=\psi_0
	\quad \mbox{a.e.\ on }\Gamma.
	\label{weak6}
\end{align}
}

We note that, due to the lack of the regularities of time derivatives, the equations 
(\ref{LW1}) and (\ref{LW5}) are replaced by the variational formulations (\ref{weak1}) and (\ref{weak4}), respectively. 
Moreover, the boundary condition (\ref{LW3}) is hidden in the weak form  (\ref{weak1}).\smallskip

Next, we define the strong solution of problem \eqref{LW1}--\eqref{LW7}. \smallskip

\paragraph{\bf Definition~2.2.} 
{\it The sextuplet $(\phi, \mu, \xi, \psi, w, \zeta)$ is called \pier{a} strong solution of problem \eqref{LW1}--\eqref{LW7} if
\begin{align*}
	& \phi \in W^{1,\infty} (0,T;V^*)\cap H^1 (0,T;V) \cap L^\infty \bigl( 0,T;H^2(\Omega) \bigr), \\
	& \mu \in L^\infty (0,T;V) \cap L^2 \bigl( 0,T;W \cap H^3(\Omega) \bigr), \quad \xi \in L^\infty (0,T;H), \\
	& \psi \in W^{1,\infty} (0,T;V_\Gamma^*)\cap H^1 (0,T;V_\Gamma) \cap L^\infty \bigl( 0,T;H^2(\Gamma) \bigr), \\
	& w \in L^\infty (0,T;V_\Gamma) \cap L^2 \bigl(0,T;H^3(\Gamma) \bigr), \quad \zeta \in L^\infty (0,T;H_\Gamma)
\end{align*} 
and they satisfy {\rm (\ref{LW1})--(\ref{LW7})}. }\smallskip

Before we state our main theorems, we recall the structure of mass conservation of problem \eqref{LW1}--\eqref{LW7}. 
\pier{Taking $z=1$ in \eqref{weak1} 
and integrating from $0$ to $t$ with the help of \eqref{weak6},}
we obtain the first equality in (\ref{mass}). 
Analogously, from \pier{(\ref{weak4}) and (\ref{weak6})} 
we obtain the second condition in  (\ref{mass}).
Therefore, it is useful to define the following mean value functions: 
\begin{align}
	& m_\Omega(z):=\frac{1}{|\Omega|}\int_{\Omega} z dx, 
	\quad |\Omega|:=\int_{\Omega }1 dx, \label{pier1}
	\\
	& m_\Gamma (z_\Gamma):=\frac{1}{|\Gamma|}\int_{\Gamma} z_\Gamma d\Gamma, 
	\quad |\Gamma |:=\int_{\Gamma }1 d\Gamma, \label{pier2}
\end{align} 
for any $z\in L^1(\Omega)$ and $z_\Gamma \in L^1(\Gamma)$. \medskip

Throughout this paper, we make the following assumptions:
\begin{enumerate}
 \item[(A1)] $\phi_0 \in V$, $\widehat{\beta}(\phi_0) \in L^1(\Omega)$, $\psi_0 \in V_\Gamma$, $\widehat{\beta}_\Gamma(\psi_0) \in L^1(\Gamma)$, and 
 ${\phi_0}_{|_\Gamma}=\psi_0$. Moreover, 
 $m_0:=m_\Omega(\phi_0) \in {\rm int}D(\beta)$, $m_{\Gamma 0}:=m_\Gamma(\psi_0) \in {\rm int}D(\beta_\Gamma)$;
 \item[(A2)] $f \in L^2(0,T;V)$, $g \in L^2(0,T;V_\Gamma)$;
 \item[(A3)] $\beta $, $\beta _{\Gamma }$ are maximal monotone graphs in 
$\mathbb{R} \times \mathbb{R}$, that is, they are the subdifferentials
\begin{gather*}
	\beta =\partial \widehat{\beta}, \quad \beta _{\Gamma }
	=\partial \widehat{\beta }_{\Gamma }
\end{gather*}
of some proper lower semicontinuous and convex functions 
$\widehat{\beta }$ and $\widehat{\beta }_{\Gamma }: \mathbb{R} \to [0,\infty ]$ 
satisfying $\widehat{\beta }(0)=\widehat{\beta}_{\Gamma }(0)=0$
with the corresponding  
effective domains denoted by $D(\beta )$ and 
$D(\beta _\Gamma)$, respectively; 
 \item[(A4)]
$\pi $, $\pi _{\Gamma }: \mathbb{R} \to \mathbb{R}$ are {L}ipschitz continuous functions 
with {L}ipschitz constants $L$ and $L_{\Gamma}$, respectively; 
 \item[(A5)] $D(\beta _\Gamma ) \subseteq D(\beta )$ and 
there exist positive constants $\varrho, c_0 >0$ such that 
\begin{gather} 
	\bigl |\beta ^\circ (r) \bigr| 
	\le \varrho \bigl |\beta _{\Gamma }^\circ (r) \bigr |+c_0
	\quad 
	\mbox{for all } r \in D(\beta _\Gamma ),
	\label{A5}
\end{gather} 
where $\beta ^\circ $ and $\beta _\Gamma ^\circ $ denote the minimal sections of 
$\beta $ and $\beta _\Gamma $.
\end{enumerate}
\smallskip

The assumption (A3) implies that 
$0 \in \beta (0)$ and $0 \in \beta _{\Gamma }(0)$. 
Moreover, the minimal section $\beta ^\circ $ of $\beta $ is defined by
$\beta ^\circ (r):=\{ r^* \in \beta (r) : |r^*|=\min _{s \in \beta (r)} |s| \}$
and same definition applies to $\beta _{\Gamma }^\circ $ in (A5). 
These assumptions are the same as in \cite{CC13, CF15}, in particular, 
the compatibility condition (A5) is essential to treat different 
potentials $\beta$ in the bulk and $\beta_\Gamma$ on the boundary. 
Of course, if one chooses $\beta=\beta_\Gamma$ to be the same potential, \pier{then} {\rm (A5)} holds automatically. \medskip

Our first result is related to the existence of global weak solutions.\smallskip

\paragraph{\bf Theorem 2.1.} 
{\it Under the assumptions {\rm (A1)}--{\rm (A5)}, 
there exists a global weak solution of problem \eqref{LW1}--\eqref{LW7} in the sense of Definition 2.1.}\smallskip

The existence of strong solutions \pier{will be} discussed in Section~4 (see Theorem 4.1). \smallskip

Our second result is the continuous dependence on the initial data and external sources, which immediately yields the uniqueness of weak solutions:\smallskip

\paragraph{\bf Theorem 2.2.} 
{\it Assume that {\rm (A3)} and {\rm (A4)} hold. Moreover, \pier{let}
$f^{(1)},f^{(2)} \in L^2(0,T;V^*)$, 
$g^{(1)},g^{(2)} \in L^2(0,T;V_\Gamma^*)$, 
$\phi_0^{(1)}, \phi_0^{(2)} \in V^*$, 
$\psi_0^{(1)}, \psi_0^{(2)} \in V_\Gamma^*$ and 
\begin{align}
	& \bigl\langle \phi_0^{(1)},1 \bigr\rangle_{V^*,V}=\bigl\langle \phi_0^{(2)},1 \bigr\rangle_{V^*,V}\pier{{}= m_0 |\Omega|}, \label{pier3}
	\\
	& \bigl\langle  \psi_0^{(1)},1 \bigr\rangle_{V^*_\Gamma,V_\Gamma}=\bigl\langle  \psi_0^{(2)},1 \bigr\rangle_{V^*_\Gamma,V_\Gamma} \pier{{}= m_{\Gamma 0} |\Gamma|}. \label{pier4}
\end{align}
Let  sextuplets of functions $(\phi^{(i)}, \mu^{(i)}, \xi^{(i)}, \psi^{(i)},w^{(i)},\zeta^{(i)})$ be weak solutions of 
problem \eqref{LW1}--\eqref{LW7} 
corresponding to the given data $f^{(i)}, g^{(i)}, \phi^{(i)}_0$ and $\psi^{(i)}_0$ for $i=1,2$. 
Then, there exists a positive constant $C>0$, depending on $L$, $L_{\Gamma}$ and $T$, such that
\begin{align*}
	& 
	\bigl| \phi^{(1)}-\phi^{(2)}
	\bigr|_{C([0,T];V^*)} 
	+
	\bigl| \psi^{(1)}-\psi^{(2)} 
	\bigr|_{C([0,T];V_\Gamma^*)} 
	+
	\bigl| \phi^{(1)}-\phi^{(2)} 
	\bigr|_{L^2(0,T;V)} 
	+
	\bigl| \psi^{(1)}-\psi^{(2)} 
	\bigr|_{L^2(0,T;V_\Gamma)} 
	\nonumber \\
	& \quad \le C \left( 
	\bigl| \phi^{(1)}_0-\phi^{(2)}_0
	\bigr|_{V^*}
	+
	\bigl| \psi^{(1)}_0-\psi^{(2)}_0
	\bigr|_{V_\Gamma^*}
	+
	\bigl| f^{(1)} - f^{(2)}
	\bigr|_{L^2(0,T;V^*)} 
	+ 
	\bigl| g^{(1)}-g^{(2)}
	\bigr|_{L^2(0,T;V_\Gamma^*)}
	\right).
\end{align*} 
}

In order to prove Theorem 2.1, \pier{we quote 
the abstract framework as in \cite{KNP95, Kub12} and} 
we also prepare the following function spaces: 
\begin{align*} 
	& V_0:=\bigl\{ z \in V : m_\Omega (z) = 0 \bigr\},  \quad 
	V_{0\ast}:=\bigl\{ z^* \in V^* : \langle z^*, 1 \rangle_{V^*,V}=0 \bigr\},   
	\\
	& V_{\Gamma,0}:=\bigl\{ z_\Gamma \in V_\Gamma : m_\Gamma (z_\Gamma) = 0 \bigr\}, 
	\quad 
	V_{\Gamma, 0\ast}:=\bigl\{ z_\Gamma^* \in V_\Gamma ^* : \langle z_\Gamma^*, 1 \rangle_{V_\Gamma^*,V_\Gamma}=0 \bigr\}.
\end{align*} 
From the {P}oincar\' e--{W}irtinger inequalities (see, e.g., \cite{Heb96}), we see that there exists a positive constant $C_{\rm P}$ such that
\begin{align*}
	& |z|_V^2 
	\le {C_{\rm P}} |z|_{V_0}^2, \quad |z|_{V_0}
	:=\left( \int_{\Omega} |\nabla z|^2 dx \right)^{1/2} 
	\quad \mbox{for all } z \in V_0, \\
	& |z_\Gamma|_{V_\Gamma}^2 
	\le {C_{\rm P}} |z_\Gamma |_{V_{\Gamma,0}}^2, \quad |z_\Gamma|_{V_{\Gamma,0}}
	:=\left( \int_{\Gamma} |\nabla_\Gamma z_\Gamma|^2 d\Gamma \right)^{1/2}
	\quad \mbox{for all }z_\Gamma \in V_{0,\Gamma}.
\end{align*} 
Then, based on the {L}ax--{M}ilgram theorem, we introduce \pier{the operator} ${\mathcal N}_\Omega: V_{0\ast} \to V_0$ by: 
$u={\mathcal N}_\Omega v$ if and only if $m_\Omega(u)=0$ and
\begin{equation}
	\int_\Omega \nabla u \cdot \nabla z dx =\langle v,z \rangle_{V^*,V} \quad \mbox{for all } z \in V;
	\label{No}
\end{equation}
analogously, \pier{we define}
${\mathcal N}_\Gamma: V_{\Gamma, 0\ast} \to V_{\Gamma,0}$ by: 
$u_\Gamma={\mathcal N}_\Gamma v_\Gamma$ if and only if $m_\Gamma(u_\Gamma)=0$ and 
\begin{equation}
	\int_\Gamma \nabla_\Gamma u_\Gamma \cdot \nabla_\Gamma z_\Gamma d\Gamma 
	=\langle v_\Gamma,z_\Gamma \rangle_{V_\Gamma^*,V_\Gamma} \quad 
	\mbox{for all } z_\Gamma \in V_\Gamma.
	\label{Ng}
\end{equation}
By virtue of these definitions, we can also introduce the norms 
\begin{equation*}
	|z|_{V_{0\ast}}:= \left( \int_\Omega | \nabla {\mathcal N}_\Omega z|^2 dx \right)^{1/2}
	\quad \mbox{for all } z \in V_{0\ast},
\end{equation*}
equivalent to the usual norm $| \cdot |_{V^*}$, for the elements of $V_{0\ast}$; and 
\begin{equation*}
	|z_\Gamma|_{V_{\Gamma,0\ast}}:=
	\left( \int_\Gamma |\nabla_\Gamma {\mathcal N}_\Gamma z_\Gamma|^2 d\Gamma \right)^{1/2}
	\quad \mbox{for all } z_\Gamma \in V_{\Gamma,0\ast},
\end{equation*}
equivalent to the usual norm $| \cdot |_{V_\Gamma^*}$, for the elements of $V_{\Gamma,0\ast}$, respectively.

%%%%% Section 3. %%%%%
\section{Well-posedness}
\setcounter{equation}{0}

In this section, we prove the existence \pier{of global weak solutions} and the continuous dependence with respect to given data. 
To do so, we introduce an approximate problem for problem \eqref{LW1}--\eqref{LW7}. The idea is based on a time-discretization scheme, the {M}oreau--{Y}osida regularization, together with a viscous {C}ahn--{H}illiard approach.

Let $N \in \mathbb{N}$ and put $h:=T/N$, the time step of discretization. 
Moreover, $\tau, \sigma \in (0,1]$ stand for  
viscosity coefficients; $\varepsilon\in(0,1]$ is used as a parameter of  {M}oreau--{Y}osida regularization for maximal monotone graphs. 
We consider the following equations and conditions for $n=0,1,\ldots,N-1$:
\begin{align}
	&\frac{\phi_{n+1}-\phi_n}{h}+\mu_{n+1}-\mu_n-\Delta \mu_{n+1} = 0 
	& \mbox{a.e.\ in } \Omega,
	\label{td1}
	\\
	&\mu_{n+1} = \tau \frac{\phi_{n+1}-\phi_n}{h} -\Delta \phi_{n+1} 
	+ \beta_\varepsilon(\phi_{n+1})+\pi(\phi_{n+1})
	- f_n	
	& \mbox{a.e.\ in } \Omega,
	\label{td2}
	\\
	&\partial _{\boldsymbol{\nu}} \mu_{n+1} = 0 
	& \mbox{a.e.\ on } \Gamma,
	\label{td3}
	\\
	& ({\phi_{n+1}})_{|_\Gamma} =\psi_{n+1} 
	& \mbox{a.e.\ on } \Gamma,
	\label{td4}
	\\
	&\frac{\psi_{n+1}-\psi_n}{h}+
	w_{n+1}-w_n - \Delta_\Gamma w_{n+1} = 0
	& \mbox{a.e.\ on } \Gamma,
	\label{td5}
	\\
	&w_{n+1} = \partial_{\boldsymbol{\nu}} \phi_{n+1} 
	+ \sigma \frac{\psi_{n+1}-\psi_n}{h}-\Delta_\Gamma \psi_{n+1} 
	+ \beta_{\Gamma, \varepsilon}(\psi_{n+1})
	& 	\nonumber\\
	&\qquad\quad  +\pi_\Gamma (\psi_{n+1}) -g_n
	& \mbox{a.e.\ on } \Gamma \pier{.}
	\label{td6}
\end{align}
Note that $\phi_0$ and $\psi_0$ are known and, in order to solve the system 
(\ref{td1})--(\ref{td6}), we need to prepare initial data $\mu_0$ and $w_0$, respectively. In the level of time-discrete approximation, we set up as follows:
\begin{equation}
	\mu_0 := 0, \quad w_0:=0. 
	\label{td7}
\end{equation}
Indeed, the terms $\mu_{n+1}-\mu_n$ in the equation (\ref{td1}) and 
$w_{n+1}-w_n$ in the equation (\ref{td5}) play a role of viscosities with the parameter $h$. 
In (\ref{td2}) and (\ref{td6}), $f_n$ and $g_n$ are known too, \pier{defined by} 
\begin{gather*}
	f_n:=\frac{1}{h} \int_{nh}^{(n+1)h} f(s)ds,
	\quad 
	g_n:=\frac{1}{h} \int_{nh}^{(n+1)h} g(s)ds \quad \mbox{for } n=0,1,\ldots,N-1.
\end{gather*}

In order to approximate the maximal monotone graphs, 
we recall the {M}oreau--{Y}osida regularization (see, e.g., \cite{Bar10, Bre73}).  
For each $\varepsilon  \in (0,1]$, we define 
$\beta _\varepsilon , \beta _{\Gamma,\varepsilon }:\mathbb{R} \to \mathbb{R}$, 
along with the associated resolvent operators 
$J_\varepsilon , J_{\Gamma,\varepsilon}:\mathbb{R} \to \mathbb{R}$ given 
by 
\begin{align*}
	& \beta _\varepsilon (r)
	:= \frac{1}{\varepsilon } \bigl( r-J_\varepsilon (r) \bigr), 
	\quad 
	J_\varepsilon (r) 
	:=(I+\varepsilon  \beta )^{-1} (r),
	\\
	& \beta _{\Gamma, \varepsilon } (r)
	:= \frac{1}{\varepsilon  \varrho} \bigl( r-J_{\Gamma,\varepsilon }(r) \bigr ), 
	\quad 
	J_{\Gamma,\varepsilon }(r):=(I+\varepsilon \varrho \beta _\Gamma )^{-1} (r),
\end{align*}
for all $r \in \mathbb{R}$, where $\varrho>0$ is same as in the condition (\ref{A5}). 
As a remark, the above two definitions are not symmetric, more precisely, 
the parameter of approximation is not directly $\varepsilon $ but $\varepsilon \varrho $ in the definition of 
$\beta_{\Gamma, \varepsilon}$ and $J_{\Gamma,\varepsilon}$. 
This is \pier{important in order to apply} \cite[Lemma~4.4]{CC13}, \pier{which ensures that} 
\begin{gather}
	\bigl |\beta_\varepsilon (r)\bigr | 
	\le \varrho \bigl |\beta _{\Gamma,\varepsilon } (r)\bigr |+c_0
	\quad 
	\mbox{for all } r \in \mathbb{R},
	\label{A5e}
\end{gather} 
for all $\varepsilon \in (0,1]$ with the same constants $\varrho $ and $c_0$ as in (\ref{A5}).  
We also have 
$\beta _\varepsilon (0)=\beta _{\Gamma, \varepsilon }(0)=0$. 
Moreover, the related {M}oreau--{Y}osida regularizations $\widehat{\beta }_\varepsilon , 
\widehat{\beta }_{\Gamma,\varepsilon }$
of $\widehat{\beta }, \widehat{\beta}_{\Gamma }:\mathbb{R} \to \mathbb{R}$ fulfill
\begin{align*}
	& \widehat{\beta }_{\varepsilon }(r)
	:=\inf_{s \in \mathbb{R}}\left\{ \frac{1}{2\varepsilon } |r-s|^2
	+\widehat{\beta }(s) \right\} 
	= 
	\frac{1}{2\varepsilon } 
	\bigl| r-J_\varepsilon  (r) \bigr|^2+\widehat{\beta }\bigl (J_\varepsilon (r) \bigr )
	= \int_{0}^{r} \beta _\varepsilon (s)ds,
	\\
	& \widehat{\beta }_{\Gamma, \varepsilon }(r)
	:=\inf_{s \in \mathbb{R}}\left\{ \frac{1}{2\varepsilon  \varrho } |r-s|^2
	+\widehat{\beta }_\Gamma (s) \right\} 
	= \int_0^r \beta _{\Gamma, \varepsilon }(s)ds,
\end{align*}
for all $r\in \mathbb{R}$. Then, we see that $\beta_\varepsilon $ and $\beta_{\Gamma, \varepsilon }$ are 
{L}ipschitz continuous with constants $1/\varepsilon$ and   
$1/(\varepsilon  \varrho)$, respectively. 
Additionally, we also use the following facts:
\begin{align}
	& \bigl |\beta _\varepsilon (r) \bigr | \le \bigl |\beta ^\circ (r) \bigr |, \quad 
	\bigl |\beta _{\Gamma ,\varepsilon }(r) \bigr | \le \bigl |\beta _{\Gamma }^\circ (r) \bigr |,
	\nonumber \\
	& 0 \le \widehat{\beta }_\varepsilon (r) \le \widehat{\beta }(r), \quad 
	0 \le \widehat{\beta }_{\Gamma,\varepsilon } (r) \le \widehat{\beta }_{\Gamma }(r),
	\label{prim}
\end{align}
for all $r\in \mathbb{R}$.

%%%%% Section 3.1. %%%%%
\subsection{Time-discrete approximate solution}

In this subsection, firstly we discuss the existence of \pier{solutions to the} time-discrete approximate problem 
(\ref{td1})--(\ref{td6}) for all $n=0,1,\ldots, N-1$, for arbitrary but fixed parameters $\tau, \sigma \in (0,1]$. 
Secondly, \pier{by introducing the piecewise linear and constant interpolants, we construct} the approximate problem of a viscous {C}ahn--{H}illiard system.\smallskip

\paragraph{\bf Proposition 3.1.} {\it \pier{There is a value $h^* \in (0,1],$ depending on $\tau$ and $\sigma$, such that for every  $h\in (0,h^*]$, 
there exists a unique quadruplet $(\phi_{n+1} ,\mu_{n+1} ,\psi_{n+1}, w_{n+1})$
with} $\phi_{n+1} \in H^2(\Omega)$, 
$\mu_{n+1} \in W$, 
$\psi_{n+1}, w_{n+1} \in H^2(\Gamma)$ such that 
{\rm (\ref{td1})}--{\rm (\ref{td6})} holds for all $n=0,1,2,\ldots, N-1$. }\smallskip

\paragraph{\bf Proof.} Define $\Delta _{\rm N}:W \to H$ be the {L}aplace operator, subject to the homogeneous {N}eumann boundary 
condition. 
From (\ref{td1}) and (\ref{td3}), we infer that
\begin{equation}
	\mu_{n+1}=(I - \Delta_{\rm N})^{-1}\left( \mu_n-\frac{\phi_{n+1}-\phi_n}{h}\right) 
	\quad \mbox{in } H,
	\label{td1b}
\end{equation}
where $I-\Delta_{\rm N}$ is a linear operator from \pier{its domain $W\subset H$} to $H$. At the same time, from (\ref{td5}) we 
obtain 
\begin{equation}
	w_{n+1}=(I - \Delta_\Gamma)^{-1}\left( w_n-\frac{\psi_{n+1}-\psi_n}{h}\right) 
	\quad \mbox{in }H_\Gamma,
	\label{td5b}
\end{equation}
where $I-\Delta_\Gamma$ is a linear operator from \pier{$H^2(\Gamma)\subset H_\Gamma$} to $H_\Gamma$. 
As a consequence, equation (\ref{td2}) can be rewritten as
\begin{align}
	& (I - \Delta_{\rm N})^{-1}\phi_{n+1} + \tau \phi_{n+1}-h \Delta \phi_{n+1} + h \beta_\varepsilon(\phi_{n+1})
	+ h \pi (\phi_{n+1}) 
	\nonumber \\
	& \quad = (I - \Delta_{\rm N})^{-1}\phi_n + h (I - \Delta_{\rm N})^{-1}\mu_n + \tau \phi_n 
	+ h f_n
	\quad \mbox{in } H
	\label{pier5}
\end{align}
and the condition (\ref{td6}) becomes 
\begin{align}
	& h \partial_{\boldsymbol{\nu}}\phi_{n+1} + 
	(I - \Delta_\Gamma)^{-1}\psi_{n+1} + \sigma \pier{\psi_{n+1}} 
	- h \Delta_\Gamma \psi_{n+1} + h \beta_{\Gamma, \varepsilon}(\psi_{n+1})
	+ h \pi_\Gamma (\psi_{n+1}) 
	\nonumber \\
	& \quad = (I - \Delta_\Gamma)^{-1}\psi_n + h (I - \Delta_\Gamma)^{-1}w_n + \sigma \psi_n 
	+ h g_n
	\quad \mbox{in } H_\Gamma.
	\label{pier6}
\end{align}

Now, the \pier{map} 
\begin{equation*}
	(z,z_\Gamma) \mapsto (-h \Delta z, h \partial _{\boldsymbol{\nu}}z - h \Delta_\Gamma z_\Gamma)
\end{equation*}
gives a maximal monotone operator $A$ from its domain $D(A):=
\{ (z,z_\Gamma) \in H^2(\Omega) \times H^2(\Gamma) : z_{|_\Gamma}=z_\Gamma \mbox{ a.e.\ on } \Gamma \}$ to 
$\boldsymbol{H}:=H \times H_\Gamma$, \pier{with} reference to 
\cite[p.~47, Theorem~2.8]{Bar10}. 
Indeed, it coincides with the subdifferential of 
\pier{the} proper, lower semicontinuous and convex functional 
$J: \boldsymbol{H} \to [0,\infty)$ defined by
\begin{equation*}
	J(z,z_\Gamma):=
	\begin{cases}
	\displaystyle 
	\frac{h}{2}\int_\Omega |\nabla z |^2 dx + \frac{h}{2} \int_\Gamma |\nabla z_\Gamma |^2 d\Gamma 
	\quad 
	\mbox{if } (z,z_\Gamma) \in \boldsymbol{V}, \\
	+\infty \quad 
	\mbox{otherwise},
	\end{cases}
\end{equation*}
where 
$\boldsymbol{V}:=\{ (z,z_\Gamma) \in V \times V_\Gamma : z_{|_\Gamma}=z_\Gamma \mbox{ a.e.\ on } \Gamma\}$. 
\pier{This} also implies that the \pier{subdifferential} of $J$ 
\pier{in $\boldsymbol{H}$} at $(z,z_\Gamma)$ coincides with 
$A(z,z_\Gamma)= (-h \Delta z, h \partial _{\boldsymbol{\nu}}z - h \Delta_\Gamma z_\Gamma)$. 
Next, we define another operator $B:\boldsymbol{H} \to \boldsymbol{H}$ by
\begin{align*}
	B(z,z_\Gamma) := & \bigl(  (I - \Delta_{\rm N})^{-1}z + \tau z+ h \beta_\varepsilon(z)
	+ h \pi (z), \\
	& \quad 
	(I - \Delta_\Gamma)^{-1}z_\Gamma+ \sigma z_\Gamma
	+ h \beta_{\Gamma, \varepsilon}(z_\Gamma)
	+ h \pi_\Gamma (z_\Gamma)  \bigr)
\end{align*}
with its domain $D(B)=\boldsymbol{H}$. Then, \pier{we see that} $B$ is {L}ipschitz continuous and 
monotone provided that $h$ is sufficiently small compared to $\tau$ and $\sigma$, namely $h \in (0,h^*]$ where
$h^*L < \tau /2$ and $h^* L_\Gamma <\sigma/2$:
\begin{align*}
	& \bigl( B(z_1,z_{\Gamma,1}) - B(z_2,z_{\Gamma,2}), (z_1 z_{\Gamma,1})-(z_2,z_{\Gamma,2})\bigr)_{\boldsymbol{H}}\\
	& \quad \ge \frac{\tau}{2} |z_1-z_2|_H^2 + \frac{\sigma}{2}|z_{\Gamma,1}-z_{\Gamma,2}|_{H_\Gamma}^2
	\quad \mbox{for all } (z_i,z_{\Gamma,i})\in \boldsymbol{H}, \ \pier{i=1,2} \takeshi{,}
\end{align*}
\pier{of course,} $B$ is also coercive. 

Hence, 
from general theory of the maximal monotone operator
\cite[pp.~35--36, Corollaries~2.1 and 2.2]{Bar10}, 
we conclude that $\mathrm{Ran}(A+B)=\boldsymbol{H}$. 
This implies that 
for sufficiently small $h \in (0,h^*]$, 
for each $\phi_n, \mu_n \in H$ and $\psi_n, w_n \in H_\Gamma$ given by the previous step, 
there exists a unique pair \pier{$(\phi_{n+1} , \psi_{n+1}) \in \boldsymbol{V} $ solving 
(\ref{pier5}) and (\ref{pier6})}, where the uniqueness is a consequence of the strict coerciveness of $B$. Next, we can recover $\mu_{n+1} \in W$ and $w_{n+1} \in H^2(\Gamma)$ from (\ref{td1b}) and 
(\ref{td5b}), respectively. 
By comparison \pier{in} the equations (\ref{td2}) and (\ref{td6}), we also deduce that 
$\phi_{n+1} \in H^2(\Omega)$ and $\psi_{n+1} \in H^2(\Gamma)$, using the elliptic regularity theory (see, e.g., \cite[Lemma A.1]{MZ05}). 
Thus, we can complete the proof of Proposition 3.1 by iterating from $n=0$ to $n=N-1$. \hfill $\Box$\medskip

\pier{According to} the standard manner, we now define the following piecewise linear functions and step functions:
\begin{align*}
	\hat\phi_h(t)& :=\phi_n + \frac{\phi_{n+1}-\phi_n}{h}(t-nh) \quad & \mbox{for } t \in \bigl[ nh, (n+1) h \bigr], 
	\quad  n=0,1,\ldots,N-1, \\
	\bar\phi_h(t)& :=\phi_{n+1} \quad & \mbox{for } t \in \bigl( nh, (n+1) h \bigr], 
	\quad  n=0,1,\ldots,N-1, \\
	\underline{f}~\! \!_h(t)& :=f_n \quad & \mbox{for } t \in \bigl( nh, (n+1) h \bigr], 
	\quad  n=0,1,\ldots,N-1, 
\end{align*}
and analogously for $\hat\mu_h, \bar\mu_n, \hat\psi_h, \bar\psi_h, \hat w_h, \bar w_h, \pier{\underline{g}~\! \!_h}$. 
Then, we have the following useful properties:
\begin{align}
	& \bigl|\hat\phi_h \bigr|_{L^2(0,T;X)}^2 \le \frac{h}{2} |\phi_0|_X^2+\bigl|\bar\phi_h \bigr|_{L^2(0,T;X)}^2, 
	\label{tool1}\\
	& \bigl|\hat\phi_h \bigr|_{L^\infty(0,T;X)} = \max\bigl\{ |\phi_0|_X, \bigl|\bar\phi_h \bigr|_{L^\infty(0,T;X)} \bigr\}, 
	\label{tool2}\\
	& \bigl|\hat\phi_h - \bar\phi_h\bigr|_{L^2(0,T;X)}^2 = \frac{h^2}{3} \bigl|\partial _t \hat \phi_h \bigr|_{L^2(0,T;X)}^2,
	\label{tool3}
\end{align}
for some suitable function space $X$. 
Indeed, (\ref{tool2}) is clear from the definition, 
\pier{the} equality (\ref{tool3}) is obtained from the direct calculation as follows:
\begin{align*}
	\bigl|\hat\phi_h -\bar\phi_h\bigr|_{L^2(0,T;X)}^2 
	& = \sum_{n=0}^{N-1} \left| \frac{\phi_{n+1}-\phi_n}{h}\right|_X^2 
	\int_{nh}^{(n+1)h} \bigl(t-h(n+1)\bigr)^2dt \\
	& = \frac{h^3}{3}\sum_{n=0}^{N-1} \left| \frac{\phi_{n+1}-\phi_n}{h}\right|_X^2  \\
	& = \frac{h^2}{3} \bigl|\partial _t \hat \phi_h \bigr|_{L^2(0,T;X)}^2.
\end{align*}
Concerning the inequality (\ref{tool1}), invoking the convexity and Jensen's inequality we obtain that 
\begin{align*}
	\bigl|\hat\phi_h \bigr|_{L^2(0,T;X)}^2 
	& = \sum_{n=0}^{N-1} \int_{nh}^{(n+1)h} 
	\left| \left( \frac{t-nh}{h} \right)\phi_{n+1} + \left[ 1- \left( \frac{t-nh}{h} \right) \right] \phi_n  \right|_X^2 dt \\
	& \le \sum_{n=0}^{N-1} \int_{nh}^{(n+1)h} \left\{ 
	\left( \frac{t-nh}{h} \right)
	| \phi_{n+1} |_X^2  + \left[ 1- \left( \frac{t-nh}{h} \right) \right] | \phi_n |_X^2 \right\} dt \\
	& = \sum_{n=0}^{N-1} \left( \frac{h}{2} | \phi_{n+1} |_X^2
	  + \frac{h}{2} | \phi_n |_X^2\right)  \\
	 & \le \frac{h}{2}|\phi_0|_X^2 +  \bigl|\bar\phi_h \bigr|_{L^2(0,T;X)}^2.
\end{align*} 
Under these settings, we see from (\ref{td1})--(\ref{td6}) that the functions 
$$
\hat\phi_h,\ \ \bar\phi_h,\ \  \hat\mu_h,\ \  \bar\mu_n,\ \  \hat\psi_h,\ \  \bar\psi_h,\ \  \hat w_h,\ \  \bar w_h
$$
constructed above solve the following polygonal approximate problem of \pier{the} viscous {C}ahn--{H}illiard system:
\begin{align} 
	&\partial_t \hat\phi_h +h \partial _t \hat \mu_h - \Delta \bar\mu_h =0 
	& \mbox{a.e.\ in }Q, 
	\label{pLW1}
	\\
	&\bar\mu_h = \tau \partial _t \hat \phi_h -\Delta \bar\phi_h 
	+ \beta_\varepsilon(\bar \phi_h) + \pi (\bar \phi_h) -\underline{f}\! {_h}
	& \mbox{a.e.\ in }Q, 
	\label{pLW2}
	\\
	&\partial_{\boldsymbol{\nu}} \bar \mu_h =0 
	& \mbox{a.e.\ on }\Sigma, 
	\label{pLW3}
	\\
	&(\bar \phi_h)_{|_{\Gamma}} =\bar \psi_h 
	& \mbox{a.e.\ on }\Sigma,
	\label{pLW4}
	\\ 
	&\partial_t \hat \psi_h + h \partial_t \hat w_h -\Delta_\Gamma \bar w_h = 0 
	& \mbox{a.e.\ on }\Sigma, 
	\label{pLW5}
	\\
	& \bar w_h = \partial _{\boldsymbol{\nu}} \bar \phi_h + \sigma \partial _t \hat \psi_h - \Delta_\Gamma \bar \psi_h 
	+ \beta_{\Gamma, \varepsilon}(\bar \psi_h) + \pi_\Gamma(\bar \psi_h)-\underline{g} {_h}
	& \mbox{a.e.\ on }\Sigma,
	\label{pLW6}
	\\
	&\hat \phi_h (0)=\phi_0, \quad \hat \mu_h (0)=0
	\quad 
	\mbox{a.e.\ in }\Omega, 
	\quad 
	\hat \psi_h (0)=\psi_0, \quad \hat w_h(0)=0
	& \mbox{a.e.\ on }\Gamma,
	\label{pLW7}
\end{align}
for every $h\in (0, h^*]$.
By virtue of the definitions of $f_n$ and $g_n$ we see that 
$\{\underline{f}\! {_h} \}_{h >0}$ and $\{\underline{g}{_h} \}_{h>0}$ are bounded in 
$L^2(0,T;V)$ and $L^2(0,T;V_\Gamma)$, respectively. Indeed, from the {H}\"older inequality \pier{we infer that}
\begin{align}
	\int_0^T \bigl| \underline{f}\! {_h}(t) \bigr|_V^2 dt 
	& = \sum_{n=0}^{N-1} h |f_n|_V^2 \nonumber
	\\
	& = \sum_{n=0}^{N-1} h \left| \frac{1}{h}\int_{nh}^{(n+1)h} f(s)ds \right|_V^2 \nonumber
	\\
	& \le \frac{1}{h} \sum_{n=0}^{N-1}  \left( \int_{nh}^{(n+1)h} \bigl| f(s) \bigr|_V^2 ds \right) \left( \int_{nh}^{(n+1)h} 1 ds \right) \nonumber
	\\
	& = |f|_{L^2(0,T;V)}^2, \quad \text{for all }h>0, 
	\label{pier7}
\end{align}
and a similar result holds for $\{\underline{g}{_h} \}_{h>0}$. 
In the next subsection we will proceed to derive necessary uniform estimates for problem (\ref{pLW1})--(\ref{pLW7}).

%%%%% Section 3.2. %%%%%
\subsection{A priori estimates and limiting procedure}

Hereafter, we derive uniform estimates that are independent of $h=T/N$ for 
problem (\ref{pLW1})--(\ref{pLW7}). 
We also take care of the dependence 
with respect to $\tau, \sigma, \varepsilon \in (0,1]$. \medskip

\paragraph{\bf Lemma 3.1.}
{\it There exists a positive constant $M_1$, 
independent of $h\in(0,h^{**}]$, $\tau, \sigma, \varepsilon \in (0,1]$, 
such that
\begin{align}
	& \bigl|\bar{\phi}_h\bigr|_{L^\infty(0,T;V)}^2+ \bigl|\bar{\psi}_h\bigr|_{L^\infty(0,T;{V_\Gamma})}^2 +
	\bigl|\partial_t \hat{\phi}_h + h \partial _t \hat{\mu}_h \bigr|_{L^2(0,T;V^*)}^2 
	+ \bigl|\partial_t \hat{\psi}_h + h \partial _t \hat{w}_h \bigr|_{L^2(0,T;{V_\Gamma}^*)}^2 
	\nonumber \\
	&\quad
	+ h | \bar{\mu}_h|_{L^\infty(0,T;H)}^2 
	+ h^2 \bigl| \partial _t \hat{\mu}_h\bigr|_{L^2(0,T;H)}^2 
	+ h | \bar{w}_h|_{L^\infty(0,T;H_\Gamma)}^2 
	+ h^2 \bigl| \partial _t \hat{w}_h\bigr|_{L^2(0,T;H_\Gamma)}^2 
	\nonumber \\
	&\quad
	+ \tau \bigl|\partial_t \hat{\phi}_h \bigr|_{L^2(0,T;H)}^2 
	+ \sigma \bigl|\partial_t \hat{\psi}_h \bigr|_{L^2(0,T;H_\Gamma)}^2 
	+ h \bigl|\partial_t \hat{\phi}_h \bigr|_{L^2(0,T;V)}^2 
	+ h \bigl|\partial_t \hat{\psi}_h \bigr|_{L^2(0,T;V_\Gamma)}^2 
	\nonumber \\
	&\quad
	+ \bigl|\pier{\widehat{\beta}_\varepsilon} (\bar{\phi}_h) \bigr|_{L^\infty(0,T;L^1(\Omega))}
	+ \bigl|\pier{\widehat{\beta}_{\Gamma, \varepsilon}} (\bar{\psi}_h) \bigr|_{L^\infty(0,T;L^1(\Gamma))}
	\le M_1,
	\label{M1}
\end{align}
for all $h\in(0,h^{**}]$, where $h^{**}\in(0,h^*]$ is a threshold \pier{value for  the step size depending} on $\tau, \sigma \in (0,1]$.} \smallskip

\paragraph{\bf Proof.} \pier{By integrating (\ref{td1}) over $\Omega$,  with the help of (\ref{td3}) we deduce the following relation for the mean values defined in \eqref{pier1}:}
\begin{equation}
	m_\Omega (\phi_{n+1}+h \mu_{n+1})=\pier{m_\Omega (\phi_{n}+h\mu_{n})} =m_0, 
	\label{m0}
\end{equation}
for all $n =0,1,\ldots,N-1$, where \pier{(cf.~{\rm (A1)} and (\ref{td7}))}
\begin{equation*}
	m_0:=m_\Omega(\phi_0)
	= \frac{1}{|\Omega|} \int_\Omega \phi_0 dx 
	= \frac{1}{|\Omega|} \int_\Omega (\phi_0+h\mu_0) dx.
\end{equation*} 
In a similar manner, from (\ref{td5}) \pier{and \eqref{pier2}} it follows that
\begin{equation}
	m_\Gamma (\psi_{n+1}+h w_{n+1})= \pier{m_\Gamma (\psi_{n}+h w_{n})}=m_{\Gamma 0}, 
	\label{mg0}
\end{equation}
for all $n =0,1,\ldots,N-1$, where, \pier{thanks to {\rm (A1)} and (\ref{td7}),}
\begin{equation*}
	m_{\Gamma 0}:=m_\Gamma(\psi_0)
	= \frac{1}{|\Gamma|} \int_\Gamma \psi_0 d\Gamma 
	= \frac{1}{|\Gamma|} \int_\Gamma (\psi_0+h w_0) d\Gamma.
\end{equation*} 
Then, \pier{as (\ref{m0}) entails} that
\begin{align*}
\langle \phi_{n+1}+h \mu_{n+1}-\phi_n-h \mu_n, 1 \rangle_{V^*,V}
& = \int_\Omega (
\phi_{n+1}+h \mu_{n+1}) dx 
- \int_\Omega (
\phi_n+h \mu_n ) dx
\pier{{} =0},
\end{align*}
that is, $\phi_{n+1}+h \mu_{n+1}-\phi_n-h \mu_n \in V_{0*}$, \pier{we can test} (\ref{td1}) 
by 
\begin{equation*}
	{\mathcal N}_\Omega (\phi_{n+1}+h \mu_{n+1}-\phi_n-h \mu_n)
\end{equation*}
\pier{and, using (\ref{td3}) and (\ref{No}),} we obtain
\begin{equation}
	h \left| \frac{\phi_{n+1}+h \mu_{n+1}-\phi_n-h \mu_n}{h}\right|_{V_{0\ast}}^2 
	+ \int_\Omega \mu_{n+1}  (\phi_{n+1}+h \mu_{n+1}-\phi_n-h \mu_n) dx = 0.
	\label{Lem31-a}
\end{equation}
for all $n =0,1,2,\ldots,N-1$. Next, 
from (\ref{mg0}) we see that
\begin{align*}
\langle \psi_{n+1}+h w_{n+1}-\psi_n-h w_n, 1 \rangle_{V_\Gamma^*,V_\Gamma}
& = \int_\Gamma (
\psi_{n+1}+h w_{n+1}) d\Gamma
- \int_\Gamma (
\psi_n+h w_n ) d\Gamma
\pier{{} =0},
\end{align*}
that is, $\psi_{n+1}+h w_{n+1}-\psi_n-h w_n \in V_{\Gamma,0*}$; therefore,
testing (\ref{td5}) 
by
\begin{equation*}
{\mathcal N}_\Gamma (\psi_{n+1}+h w_{n+1}-\psi_n-h w_n),
\end{equation*}
 and using (\ref{Ng}) we obtain
\begin{equation}
	h \left| \frac{\psi_{n+1}+h w_{n+1}-\psi_n-h w_n}{h}\right|_{V_{\Gamma, 0\ast}}^2 
	+ \int_\Gamma w_{n+1}  (\psi_{n+1}+h w_{n+1}-\psi_n-h w_n) dx = 0,
	\label{Lem31-b}
\end{equation}
for all $n =0,1,\ldots,N-1$. 
Next, we add $\phi_{n+1}$ to both sides of (\ref{td2}), 
multiply the resultant by $\phi_{n+1}-\phi_n$ and use the condition (\ref{td4}) and the equation
(\ref{td6}), to find out that
\begin{align}
	& \int_\Omega \mu_{n+1}(\phi_{n+1}-\phi_n) dx 
	+ \int_\Gamma w_{n+1}(\psi_{n+1}-\psi_n) d\Gamma
	\nonumber \\
	& \ge  \tau h
	\left| \frac{\phi_{n+1}-\phi_n}{h}\right|_H^2
	+  \sigma h \left| \frac{\psi_{n+1}-\psi_n}{h}\right|_{H_\Gamma}^2
	+ \frac{1}{2}|\phi_{n+1}|^2_V+\frac{1}{2}|\phi_{n+1}-\phi_n|_V^2- \frac{1}{2}|\phi_n|^2_V 
	\nonumber 	
	\\
	& 
	\quad {}+ \frac{1}{2}|\psi_{n+1}|^2_{V_\Gamma}+\frac{1}{2}|\psi_{n+1}-\psi_n|_{V_\Gamma}^2- \frac{1}{2}|\psi_n|^2_{V_\Gamma}
	+ \int_\Omega \widehat\beta_\varepsilon(\phi_{n+1}) dx 
	-  \int_\Omega \widehat\beta_\varepsilon(\phi_n) dx
	\nonumber 
	\\
	& 
	\quad {} + \int_\Gamma \widehat\beta_{\Gamma,\varepsilon}(\psi_{n+1}) d\Gamma 
	-  \int_\Gamma \widehat\beta_{\Gamma,\varepsilon}(\psi_n) d\Gamma
	+ \int_\Omega \bigl( \pi(\phi_{n+1}) - f_n - \phi_{n+1} \bigr)(\phi_{n+1}-\phi_n)dx 
	\nonumber 
	\\
	& \quad {} 
	+ \int_\Gamma \bigl( \pi_\Gamma (\psi_{n+1}) - g_n - \psi_{n+1} \bigr)(\psi_{n+1}-\psi_n)d\Gamma 
	\label{Lem31-c}
\end{align} 
for all $n =0,1,\ldots,N-1$, 
where we used \pier{the elementary inequality} $r(r-s)=(r^2+(r-s)^2-s^2)/2$ for $r,s\in\mathbb{R}$. 

Now, we collect (\ref{Lem31-a})--(\ref{Lem31-c}), sum up for 
$n=0,1,\ldots,m-1$ and apply (\ref{td7}) and (\ref{prim}), obtaining
\begin{align}
	& 
	\sum_{n=0}^{m-1} h 
	\left| \frac{\phi_{n+1}-\phi_n}{h}+\mu_{n+1}- \mu_n\right|_{V_{0\ast}}^2 
	+
	\sum_{n=0}^{m-1} h 
	\left| \frac{\psi_{n+1}-\psi_n}{h}+w_{n+1}- w_n\right|_{V_{\Gamma,0\ast}}^2
	\nonumber \\
	& \quad {}
	+ 
	\frac{h}{2} |\mu_m|_H^2 
	+
	\sum_{n=0}^{m-1} \frac{h}{2}|\mu_{n+1}-\mu_n|_H^2 
	+ 
	\frac{h}{2} |w_m|_{H_\Gamma}^2 
	+ 
	\sum_{n=0}^{m-1} \frac{h}{2}|w_{n+1}-w_n|_{H_\Gamma}^2
	\nonumber \\
	& \quad {}
	+ \tau \sum_{n=0}^{m-1}\pier{h} \left| \frac{\phi_{n+1}-\phi_n}{h}\right|_H^2
	+ \sigma \sum_{n=0}^{m-1}\pier{h} \left| \frac{\psi_{n+1}-\psi_n}{h}\right|_{H_\Gamma}^2
	+ \frac{1}{2}|\phi_m|^2_V+\frac{1}{2} \sum_{n=0}^{m-1} |\phi_{n+1}-\phi_n|_V^2
	\nonumber 
	\\
	& 
	\quad {}+ \frac{1}{2}|\psi_m|^2_{V_\Gamma}
	+\frac{1}{2}\sum_{n=0}^{m-1} 
	|\psi_{n+1}-\psi_n|_{V_\Gamma}^2
	+ \int_\Omega \widehat\beta_\varepsilon(\phi_m) dx 
	+ \int_\Gamma \widehat\beta_{\Gamma,\varepsilon}(\psi_m) d\Gamma 
	\nonumber \\
	& \le  \frac{1}{2}|\phi_0|_V^2+\frac{1}{2}|\psi_0|_{V_\Gamma}^2
	+ \int_\Omega \widehat\beta(\phi_0) dx 
	+ \int_\Gamma \widehat\beta_{\Gamma}(\psi_0) d\Gamma 
	\nonumber \\
	& \quad {}
	\pier{-} \sum_{n=0}^{m-1} \int_\Omega \bigl( \pi(\phi_{n+1}) - f_n - \phi_{n+1} \bigr)(\phi_{n+1}-\phi_n)dx 
	\nonumber 
	\\
	& \quad {} \pier{-} \sum_{n=0}^{m-1}
	\int_\Gamma \bigl( \pi_\Gamma (\psi_{n+1}) - g_n - \psi_{n+1} \bigr)(\psi_{n+1}-\psi_n)d\Gamma 
	\label{Lem31-d}
\end{align} 
for all $m=1,2,\ldots,N$. 
We \pier{know} that there exists a positive constant $C_1$ such that 
$|z|_{V^*}^2 \le C_1 |z|_{V_{0\ast}}^2$
for all $z \in V_{0\ast}$\pier{, as well as  $|z|_{V^*}^2 \le C_1 |z|_{H}^2$ for all $z \in H$.}
Therefore, in \pier{order to estimate the terms on the right hand side of \eqref{Lem31-d} we can argue with the
help of assumptions {\rm (A2)} and {\rm (A4)}. First, we have that}
\begin{align}
	& \left| \sum_{n=0}^{m-1} \int_\Omega \bigl( \pi(\phi_{n+1}) - f_n - \phi_{n+1} \bigr)(\phi_{n+1}-\phi_n)dx \right|
	\nonumber \\
	& \le \delta \sum_{n=0}^{m-1} h \left| \frac{\phi_{n+1}-\phi_n}{h}\right|_{V^*} ^2 
	+C_\delta \sum_{n=0}^{m-1} h \bigl( 1+ |\phi_{n+1}|_V^2+ |f_n|_V^2 \bigr)
	\nonumber \\ 
	& \le 2\delta C_1 \sum_{n=0}^{m-1} h \left| \frac{\phi_{n+1}-\phi_n}{h}+\mu_{n+1}-\mu_n \right|_{V_{0\ast}} ^2 
	+2 \delta \pier{C_1} \sum_{n=0}^{m-1} h \left| \mu_{n+1}-\mu_n \right|_H ^2
	\nonumber \\
	& \quad {}+C_\delta \sum_{n=0}^{m-1} h \bigl( 1+ |\phi_{n+1}|_V^2+ |f_n|_V^2 \bigr),
	\label{Lem31-e}
\end{align}
for all $\delta>0$, 
where we also use {Y}oung's inequality with $\delta>0$; $C_\delta$ is a positive constant \pier{such} that 
$C_\delta \to \infty$ as $\delta \to 0$. 
Indeed, taking care of (\ref{m0}), we have $(\phi_{n+1}-\phi_n)/h+\mu_{n+1}-\mu_n \in V_{0\ast}$ for all 
$n=0,1,\ldots,N-1$. From (\ref{mg0}), a very similar procedure can be used to estimate the other contribution
\begin{align}
	& \left| \sum_{n=0}^{m-1} \int_\Gamma 
	\bigl( \pi_\Gamma(\psi_{n+1}) - g_n - \psi_{n+1} \bigr)(\psi_{n+1}-\psi_n)d\Gamma \right|
	\nonumber \\
	& \le 2 \delta C_2
	\sum_{n=0}^{m-1} h \left| \frac{\psi_{n+1}-\psi_n}{h}+w_{n+1}-w_n \right|_{V_{\Gamma, 0\ast}} ^2 
	+2 \delta \pier{C_2}  \sum_{n=0}^{m-1} h \left| w_{n+1}-w_n \right|_{H_\Gamma} ^2
	\nonumber \\
	& \quad {}+C_\delta \sum_{n=0}^{m-1} h \bigl( 1+ |\psi_{n+1}|_{V_\Gamma}^2+ |g_n|_{V_\Gamma}^2 \bigr),
	\label{Lem31-f}
\end{align} 
where $C_2$ is a positive constant such that 
$|z_\Gamma|_{V_\Gamma^*}^2 \le C_2 |z_\Gamma|_{V_{\Gamma,0\ast}}^2$
for all $z_\Gamma \in V_{\Gamma,0\ast}$
\pier{and $|z_\Gamma|_{V_\Gamma^*}^2 \le C_2 |z_\Gamma|_{H_{\Gamma}}^2$
for all $z_\Gamma \in H_{\Gamma}$}.
Then, we can choose $\delta>0$ in order \pier{that} 
$$
\delta\le \min\pier{\{1/(8C_1),1/(8C_2)\}}
$$
and consequently we \pier{fix} the constant $C_\delta$ in the above estimates. 
Next, we choose 
a threshold \pier{value for} the step size $h^{**}\in (0,h^*]$ \pier{with the requirement that}
\begin{equation*}
	C_\delta \pier{h^{**}} \le \frac{1}{4}. 
\end{equation*}
Then, collecting \pier{(\ref{Lem31-d})--(\ref{Lem31-f}) and recalling \eqref{pier3}--\eqref{pier4}, it is not difficult to} obtain 
\begin{align*} 
	& 
	\frac{1}{2}|\phi_m|^2_V+ \frac{1}{2}|\psi_m|^2_{V_\Gamma}
	\nonumber \\
	&\quad  \le  \frac{1}{2}|\phi_0|_V^2+\frac{1}{2}|\psi_0|_{V_\Gamma}^2
	+ \int_\Omega \widehat\beta(\phi_0) dx 
	+ \int_\Gamma \widehat\beta_{\Gamma}(\psi_0) d\Gamma 
	\nonumber \\
	&\qquad 
	+ C_\delta \sum_{n=0}^{m-1} h \bigl( 1+  |f_n|_V^2 \bigr)
	+ C_\delta \sum_{n=1}^{m-1} h |\phi_n|_V^2
	+ C_\delta h |\phi_m|_V^2 \nonumber\\
	&\qquad
	+ C_\delta \sum_{n=0}^{m-1} h \bigl( 1+ |g_n|_{V_\Gamma}^2 \bigr) 
	+ C_\delta \sum_{n=1}^{m-1} h |\psi_n|_{V_\Gamma}^2
	+ C_\delta h |\psi_m|_{V_\Gamma}^2
\end{align*}
\pier{and consequently}
\begin{align*}
	& 
	\frac{1}{2}|\phi_m|^2_V+ \frac{1}{2}|\psi_m|^2_{V_\Gamma}
	\nonumber \\
	&\quad 
	 \le  \frac{1}{2}|\phi_0|_V^2+\frac{1}{2}|\psi_0|_{V_\Gamma}^2
	+ \int_\Omega \widehat\beta(\phi_0) dx 
	+ \int_\Gamma \widehat\beta_{\Gamma}(\psi_0) d\Gamma \nonumber\\
	&\qquad 
	+ C_\delta \left( T + |f|_{L^2(0,T;V)}^2 \right)
	+ C_\delta \sum_{n=0}^{m-1} h |\phi_n|_V^2
	+ \frac{1}{4}|\phi_m|_V^2 \nonumber \\
	& \qquad {}
	+ C_\delta \left(T + |g|_{L^2(0,T;V_\Gamma)}^2 \right)
	+ C_\delta \sum_{n=0}^{m-1} h |\psi_n|_{V_\Gamma}^2
	+ \frac{1}{4}|\psi_m|_{V_\Gamma}^2,
\end{align*} 
for all $m=1,2,\ldots, N$. 
Therefore, applying the discrete {G}ronwall lemma with assumptions (A1) and (A2), we 
conclude that there exists a positive constant $\tilde{M}_1$,  independent of 
 $h\in(0,h^{**}]$, $\tau, \sigma, \varepsilon \in (0,1]$, such that
\begin{equation*}
	|\phi_m|^2_V+ |\psi_m|^2_{V_\Gamma} \le \tilde{M}_1
\end{equation*}
for all $m=1,2,\ldots, N$. 
Moreover, going back to (\ref{Lem31-d})--(\ref{Lem31-f}), we \pier{plainly deduce} 
(\ref{M1}) for some positive constant $M_1 \ge \tilde{M}_1$ independent of 
 $h\in(0,h^{**}]$ \pier{and}  $\tau, \sigma, \varepsilon \in (0,1]$. 
\hfill $\Box$\medskip

\paragraph{\bf Lemma 3.2.}  
{\it There exist two functions $\Lambda_1, \Lambda_2 \in L^2(0,T)$ and 
a positive constant $M_2$, independent of 
$h\in(0,h^{**}]$, $\tau, \sigma, \varepsilon \in (0,1]$, such that
\begin{align}
	 & \bigl| \bar{\mu}_h(t)-m_\Omega \bigl( \bar{\mu}_h(t) \bigr) \bigr|_V 
	+ 
	\bigl| \bar{w}_h(t)-m_\Gamma \bigl( \bar{w}_h(t) \bigr) \bigr|_{V_\Gamma} 
	 \le \Lambda_1(t), 
	 \label{L1}
	 \\
	 & \bigl| \beta_\varepsilon\bigl( \bar{\phi}_h(t) \bigr) \bigr|_{L^1(\Omega)} 
	+ 
	\bigl| \beta_{\Gamma, \varepsilon}\bigl( \bar{\psi}_h (t) \bigr) \bigr|_{L^1(\Gamma)} 
	 \le M_2 \left( \Lambda_2(t) + \bigl| \partial_{\boldsymbol{\nu}} \bar \phi_h (t) \bigr|_{H_\Gamma} \right)
	 \label{L2}
\end{align}
for a.a.\ $t \in (0,T)$.}\smallskip

\paragraph{\bf Proof.} 
Firstly, multiplying (\ref{pLW1}) by $\bar{\mu}_h-m_\Omega(\bar{\mu}_h)$, 
integrating the resultant over $\Omega$, using the boundary condition
(\ref{pLW3}), and applying the {Y}oung and {P}oincar\' e--{W}irtinger inequalities, we obtain 
\begin{align*}
	\int_\Omega \bigl| \nabla \bigl( \bar{\mu}_h-m_\Omega(\bar{\mu}_h) \bigr) \bigr|^2 dx 
	& \le \bigl|\partial_t \hat{\phi}_h + h \partial _t \hat{\mu}_h \bigr|_{V^*} 
	\bigl| \bar{\mu}_h-m_\Omega(\bar{\mu}_h)  \bigr|_V 
	\nonumber \\
	& 
	\le \frac{1}{2\delta}\bigl|\partial_t \hat{\phi}_h + h \partial _t \hat{\mu}_h \bigr|_{V^*} ^2 
	+ \frac{\delta}{2} \bigl| \bar{\mu}_h-m_\Omega(\bar{\mu}_h)  \bigr|_V^2 
	\nonumber \\
	& \le \frac{1}{2\delta}\bigl|\partial_t \hat{\phi}_h + h \partial _t \hat{\mu}_h \bigr|_{V^*} ^2 
	+ \frac{C_{\rm P}\delta}{2} \bigl| \bar{\mu}_h-m_\Omega(\bar{\mu}_h)  \bigr|_{V_0}^2
\end{align*}
a.e.\ in $(0,T)$\pier{, for all $\delta >0$.}
Similarly, multiplying (\ref{pLW5}) by $\bar{w}_h-m_\Gamma(\bar{w}_h)$ 
and integrating the resultant over $\Gamma$, we obtain 
\begin{align*}
	\int_\Gamma \bigl| \nabla_\Gamma \bigl( \bar{w}_h-m_\Gamma(\bar{w}_h) \bigr) \bigr|^2 d\Gamma
	&	
	\le \bigl|\partial_t \hat{\psi}_h + h \partial _t \hat{w}_h \bigr|_{V_\Gamma^*} 
	\bigl| \bar{w}_h-m_\Gamma(\bar{w}_h)  \bigr|_{V_\Gamma}
	\nonumber \\
	& \le 
	\frac{1}{2\delta}\bigl|\partial_t \hat{\psi}_h + h \partial _t \hat{w}_h \bigr|_{V_\Gamma^*} ^2 
	+ \frac{C_{\rm P}\delta}{2} \bigl| \bar{w}_h-m_\Gamma(\bar{w}_h)  \bigr|_{V_{\Gamma,0}}^2
\end{align*}
a.e.\ in $(0,T)$.
\pier{Letting} $\delta:=1/C_{\rm P}$ and applying the {P}oincar\' e--{W}irtinger inequality again, we 
deduce that 
\begin{align*}
	& \bigl| \bar{\mu}_h-m_\Omega(\bar{\mu}_h) \bigr|_V^2 
	\le 
	C_{\rm P}
	\bigl|\bar{\mu}_h-m_\Omega(\bar{\mu}_h) \bigr|_{V_0}^2 
	\le C_{\rm P}^2 
	\bigl|\partial_t \hat{\phi}_h + h \partial _t \hat{\mu}_h \bigr|_{V^*} ^2,
	\nonumber \\
	& \bigl| \bar{w}_h-m_\Gamma(\bar{w}_h) \bigr|_{V_\Gamma}^2 
	\le C_{\rm P}^2 
	\bigl|\partial_t \hat{\psi}_h + h \partial _t \hat{w}_h \bigr|_{V_\Gamma^*} ^2,
\end{align*}
a.e.\ in $(0,T)$. Taking \pier{(\ref{M1}) into account}, we conclude the estimate (\ref{L1}) with 
\begin{equation}
	\Lambda_1(t):=\sqrt{2} C_{\rm P} \left( 
	\bigl|\partial_t \hat{\phi}_h(t) + h \partial _t \hat{\mu}_h(t) \bigr|_{V^*}
	+ \bigl|\partial_t \hat{\psi}_h(t) + h \partial _t \hat{w}_h(t) \bigr|_{V_\Gamma^*}\right),
	\label{Lam1}
\end{equation}
for a.a.\ $t \in (0,T)$.

Secondly, recalling (\ref{m0}) and (\ref{mg0}), we have
\begin{gather*}
	m_\Omega \bigl(\bar\phi_h+h \bar\mu_h \bigr)=m_0,  
	\quad 
	m_\Gamma \bigl(\bar\psi_h+h \bar w_h \bigr)=m_{\Gamma 0}
	\quad \mbox{a.e.\ on } (0,T).
\end{gather*}
Then, we multiply (\ref{pLW2}) by $\bar\phi_h-m_0$ \pier{and}
integrate the resultant over $\Omega$\pier{. Also, we use (\ref{pLW4}) and exploit the argument devised in \cite[Appendix, Prop. A.1]{MZ04} (see also \cite{GMS09} for a complete proof) along} with (A1) to infer that 
there exist two positive constants $C_3, C_4>0$, independent of $h\in(0,h^{**}]$, $\tau, \sigma, \varepsilon \in (0,1]$, such that
\begin{align}
	& C_3 \int_\Omega \bigl| \beta_\varepsilon \bigl(\bar\phi_h \bigr) \bigr| dx
	-C_4
	\nonumber \\
	&  \le \int_\Omega \beta_\varepsilon \bigl( \bar\phi_h \bigr) \bigl( \bar\phi_h -m_0\bigr)dx
	\nonumber \\
	& = - \int_\Omega \bigl( \tau \partial _t \hat \phi _h + \pi \bigl( \bar \phi_h \bigr) - \underline{f}\! {_h} \bigr)
	\bigl( \bar \phi_h - m_0 \bigr) dx - \int_\Omega \bigl| \nabla \bar \phi _h \bigr|^2 dx \nonumber\\
	&\quad 
	+ \int_\Omega \bar \mu _h \bigl( \bar \phi_h + h \bar \mu_h - m_0 \bigr) dx 
	 - h \int_\Omega | \bar \mu_h |^2 dx 
	+ \int_\Gamma \partial _{\boldsymbol{\nu}} \bar \phi_h 
	\bigl( \bar \psi_h -m_{\Gamma0} \bigr) d\Gamma \nonumber 
	\\	
	& \quad {}
	+ \int_\Gamma \partial _{\boldsymbol{\nu}} \bar\phi_h (m_{\Gamma0}-m_0) d\Gamma,
	\label{Lem32-a}
\end{align}
a.e.\ in $(0,T)$. 
\pier{Similarly}, we multiply (\ref{pLW6}) by $\bar\psi_h-m_{0\Gamma}$
and integrate the resultant over $\Gamma$ to infer that 
there exist two positive constants $C_5, C_6>0$, independent of $h\in(0,h^{**}]$, $\tau, \sigma, \varepsilon \in (0,1]$, such that
\begin{align}
	& C_5 \int_\Gamma \bigl| \beta_{\Gamma,\varepsilon} \bigl(\bar\psi_h \bigr) \bigr| d\Gamma
	-C_6
	\nonumber \\
	& \le \int_\Gamma \beta_{\Gamma,\varepsilon} 
	\bigl( \bar\psi_h \bigr) \bigl( \bar\psi_h -m_{\Gamma0}\bigr)d\Gamma
	\nonumber \\
	&  = - \int_\Gamma \partial_{\boldsymbol{\nu}} 
	\bar\phi_h\bigl( \bar \psi_h - m_{\Gamma0} \bigr) d\Gamma
	 - \int_\Gamma \bigl( \sigma \partial _t \hat \psi _h 
	 + \pi_\Gamma \bigl( \bar \psi_h \bigr) - \underline{g} {_h} \bigr)
	\bigl( \bar \psi_h - m_{\Gamma0} \bigr) d\Gamma 
	\nonumber 
	\\
	& 
	\quad {}
	- \int_\Gamma \bigl| \nabla_\Gamma \bar \psi _h \bigr|^2 d\Gamma 
	+ \int_\Gamma \bar w_h \bigl( \bar \psi_h + h \bar w_h - m_{\Gamma0} \bigr) d\Gamma 
	- h \int_\Gamma | \bar w_h |^2 d\Gamma,
	\label{Lem32-b}
\end{align}
a.e.\ in $(0,T)$. 
Adding (\ref{Lem32-a}) and (\ref{Lem32-b}), we obtain that 
\begin{align}
	& C_3 \bigl| \beta_\varepsilon \bigl(\bar\phi_h \bigr) \bigr|_{L^1(\Omega)}
	+ 
	C_5 \bigl| \beta_{\Gamma,\varepsilon} \bigl(\bar\psi_h \bigr) \bigr|_{L^1(\Gamma)}
	\nonumber \\
	& \quad \le 
	C_4+C_6
	+\bigl| \tau \partial _t \hat \phi _h + \pi \bigl( \bar \phi_h \bigr) - \underline{f}\! {_h}\bigr|_H
	\bigl| \bar \phi_h - m_0 \bigr|_H 
	+ \bigl| \bar \mu _h -m_\Omega (\bar \mu_h) \bigr|_H 
	\bigl| \bar \phi_h + h \bar \mu_h - m_0 \bigr|_H 
	\nonumber 
	\\	
	& \qquad {}
	+ \bigl| \partial _{\boldsymbol{\nu}} \bar \phi_h \bigr|_{H_\Gamma}  
	| m_{\Gamma0}-m_0 |_{H_\Gamma}
	+ \bigl| \sigma \partial _t \hat \psi _h + \pi_\Gamma \bigl( \bar \psi_h \bigr) - \underline{g} {_h} \bigr|_{H_\Gamma}
	\bigl| \bar \psi_h - m_{\Gamma0} \bigr|_{H_\Gamma} 
	\nonumber \\
	& \qquad {}
	+ \bigl| \bar w_h -m_\Gamma (\bar w_h) \bigr|_{H_\Gamma} 
	\bigl| \bar \psi_h + h \bar w_h - m_{\Gamma0} \bigr|_{H_\Gamma},
	\label{Lem32-c}
\end{align}
a.e.\ in $(0,T)$. Hence, from Lemma~3.1 \pier{and \eqref{L1} it follows that} there exist a function $\Lambda_2 \in L^2(0,T)$ \pier{and}
a positive constant $M_2$, independent of 
$h\in(0,h^{**}]$, $\tau, \sigma, \varepsilon \in (0,1]$, such that
(\ref{L2}) holds. \hfill $\Box$\medskip

\paragraph{\bf Lemma 3.3.}  
{\it There exist functions $\Lambda_3, \Lambda_4 \in L^2(0,T)$ and 
positive constants $M_3, M_4$, independent of 
$h\in(0,h^{**}]$, $\tau, \sigma, \varepsilon \in (0,1]$, such that
\begin{align}
	& \bigl| m_\Omega \bigl( \bar{\mu}_h(t) \bigr) \bigr|
	+ 
	\bigl| m_\Gamma \bigl( \bar{w}_h(t) \bigr)  \bigr|
	 \le M_3 \left( \Lambda_3(t) + \bigl| \partial_{\boldsymbol{\nu}} \bar \phi_h (t) \bigr|_{H_\Gamma} \right),
	 \label{L3}
	 \\
	& \bigl| \bar{\mu}_h(t) \bigr|_V 
	+ 
	\bigl| \bar{w}_h(t) \bigr|_{V_\Gamma} 
	 \le M_4 \left( \Lambda_4(t) + \bigl| \partial_{\boldsymbol{\nu}} \bar \phi_h (t) \bigr|_{H_\Gamma} \right),
	 \label{L4}
\end{align}
for a.a.\ $t \in (0,T)$.}\smallskip

\paragraph{\bf Proof.} 
Multiplying the equation (\ref{pLW2}) 
by $1/|\Omega|$ and
the equation (\ref{pLW6}) 
by $1/|\Gamma|$\pier{,} using integration by parts, and adding the resultants together, we obtain 
\begin{align*}
	& \bigl| m_\Omega (\bar \mu_h) \bigr| 
	+ \bigl| m_\Gamma (\bar w _h ) \bigr|
	\nonumber \\
	& \le 
	\frac{1}{|\Omega|}\int_\Omega \bigl| \beta_\varepsilon(\bar \phi_h) 
	\bigr| dx 
	+
	\frac{1}{|\Omega|}\int_\Omega \bigl| \tau \partial _t \hat \phi_h 
	+  \pi (\bar \phi_h) -\underline{f}\! {_h} 
	\bigr| dx 
	+ \frac{1}{|\Omega|} \int_\Gamma \bigl| \partial_{\boldsymbol{\nu}} \bar \phi_h\bigr| d\Gamma
	\nonumber 
	\\
	& \quad {} 
	+\frac{1}{|\Gamma|} \int_\Gamma 
	\bigl| \beta_{\Gamma, \varepsilon}(\bar \psi_h) \bigr| d\Gamma
	+\frac{1}{|\Gamma|} \int_\Gamma 
	\bigl| \partial _{\boldsymbol{\nu}} \bar \phi_h + \sigma \partial _t \hat \psi_h  
	+ \pi_\Gamma(\bar \psi_h)-\underline{g} {_h} \bigr| d\Gamma,
\end{align*}
a.e.\ in $(0,T)$. Therefore, using {H}\"older's inequality and 
recalling (\ref{M1}) and (\ref{L2}), 
we see 
that there exist a function $\Lambda_3 \in L^2(0,T)$ \pier{depending} on $M_2, \Lambda_2$ and 
a positive constant $M_3$  
\pier{depending} on $M_2, L, L_\Gamma, |f|_{L^2(0,T;H)}, |g|_{L^2(0,T;H_\Gamma)}, |\Omega|$ and $|\Gamma|$,
independent of 
$h\in(0,h^{**}]$, $\tau, \sigma, \varepsilon \in (0,1]$, such that the estimate 
(\ref{L3}) holds. 
Additionally, from 
(\ref{L1}) and the above estimate on mean values of $\bar \mu_h$, $\bar w _h$,   we \pier{infer (\ref{L4}), where the function $\Lambda_4 \in L^2(0,T)$ and the
positive constant} $M_4$ are independent of 
$h\in(0,h^{**}]$, $\tau, \sigma, \varepsilon \in (0,1]$. 
\hfill $\Box$\medskip

\paragraph{\bf Lemma 3.4.}  
{\it There exist a function $\Lambda_5 \in L^2(0,T)$ and 
a positive constant $M_5$, independent of 
$h\in(0,h^{**}]$, $\tau, \sigma, \varepsilon \in (0,1]$, such that
\begin{equation}
	\bigl| \beta_\varepsilon \bigl( \bar \phi_h (t)\bigr) \bigr|_H
	+
	\bigl| \beta_{\Gamma, \varepsilon} \bigl( \bar \psi_h (t)\bigr) \bigr|_{H_\Gamma}
	\le 
	M_5 \left( \Lambda_5(t) + \bigl| \partial _{\boldsymbol{\nu}} \bar \phi _h(t) \bigr|_{H_\Gamma} \right)
	\label{L5}
\end{equation} 
for a.a.\ $t \in (0,T)$.}\smallskip

\paragraph{\bf Proof.} 
Multiplying the equation (\ref{pLW2}) 
by $\beta_\varepsilon (\bar \phi_h)$, integrating over 
$\Omega$, and using (\ref{pLW4}), we have 
\begin{align}
	& \int _\Omega \beta_\varepsilon' \bigl( \bar \phi _h \bigr) 
	\bigl| \nabla \bar \phi_h \bigr|^2 dx 
	+ 
	\int_\Omega \bigl| \beta_\varepsilon \bigl( \bar \phi_h \bigr) \bigr|^2 dx 
	\nonumber \\
	& \le 
	\frac{1}{2}\bigl| \bar\mu_h - \tau \partial_t \hat \phi _h - \pi \bigl( \bar \phi_h \bigr) 
	+\underline{f}\! {_h} \bigr|_H ^2
	+ \frac{1}{2}
	\bigl| \beta_\varepsilon \bigl( \bar \phi_h \bigr) \bigr|_H^2
	+ \int_\Gamma \partial _{\boldsymbol{\nu}} \bar \phi_h \beta_\varepsilon \bigl( \bar \psi_h \bigr)
	d\Gamma, 
	\label{Lem34-a}
\end{align} 
a.e.\ in $(0,T)$. Next, multiplying the the equation (\ref{pLW6}) 
by $\beta_{\Gamma, \varepsilon} (\bar \psi_h)$ and integrating the resultant over $\Gamma$, we have 
\begin{align}
	& \int _\Gamma \beta_{\Gamma,\varepsilon}' \bigl( \bar \psi _h \bigr) 
	\bigl| \nabla_\Gamma \bar \psi_h \bigr|^2 d\Gamma 
	+ 
	\int_\Gamma \bigl| \beta_{\Gamma, \varepsilon} \bigl( \bar \psi_h \bigr) \bigr|^2 d\Gamma 
	\nonumber \\
	&  \le 
	\bigl| \bar w_h - \sigma \partial_t \hat \psi _h - \pi_\Gamma \bigl( \bar \psi_h \bigr) 
	+\underline{g} {_h} \bigr|_{H_\Gamma} ^2
	+ \frac{1}{4}
	\bigl| \beta_{\Gamma, \varepsilon} \bigl( \bar \psi_h \bigr) \bigr|_{H_\Gamma}^2
	- \int_\Gamma \partial _{\boldsymbol{\nu}} \bar \phi_h \beta_{\Gamma, \varepsilon} \bigl( \bar \psi_h \bigr)
	d\Gamma,
	\label{Lem34-b}
\end{align} 
a.e.\ in $(0,T)$. 
Now, recalling (\ref{A5e}), we can find a positive constant 
$\tilde{M}_5$ that
\begin{align*}
	& \left| \int_\Gamma \partial _{\boldsymbol{\nu}} \bar \phi_h \beta_\varepsilon \bigl( \bar \psi_h \bigr)
	d\Gamma
	- \int_\Gamma \partial _{\boldsymbol{\nu}} \bar \phi_h \beta_{\Gamma, \varepsilon} \bigl( \bar \psi_h \bigr)
	d\Gamma \right| 
	\nonumber \\
	&  \le 
	\bigl| \partial_{\boldsymbol{\nu}} \bar \phi _h \bigr|_{H_\Gamma}
	\left( \varrho \bigl| \beta_{\Gamma,\varepsilon} \bigl( \bar \psi _h \bigr) \bigr|_{H_\Gamma} +\pier{|\Gamma|^{1/2}}c_0 + 
	\bigl| \beta_{\Gamma,\varepsilon} \bigl( \bar \psi _h \bigr) \bigr|_{H_\Gamma} \right)
	\nonumber \\
	& \le \frac{1}{4}\bigl| \beta_{\Gamma,\varepsilon} \bigl( \bar \psi _h \bigr) \bigr|_{H_\Gamma}^2
	+ \tilde{M}_5 \left( 1 + \bigl| \partial _{\boldsymbol{\nu}} \bar \phi _h \bigr|_{H_\Gamma}^2 \right),
\end{align*} 
by {Y}oung's inequality. 
Finally, adding (\ref{Lem34-a}), (\ref{Lem34-b}), 
recalling the monotonicity of $\beta_\varepsilon$ in (\ref{Lem34-a})
and $\beta_{\Gamma,\varepsilon}$ in (\ref{Lem34-b}), and 
using the above estimate, 
we 
infer that there exist \pier{a function $\Lambda_5 \in L^2(0,T)$ 
and a positive constant $M_5$, both} independent of 
$h\in(0,h^{**}]$, $\tau, \sigma, \varepsilon \in (0,1]$ such that (\ref{L5}) holds a.e.\ in $(0,T)$. \hfill $\Box$\medskip

\paragraph{\bf Lemma 3.5.}  
{\it There exists a function $\Lambda_6 \in L^2(0,T)$, independent of 
$h\in(0,h^{**}]$, $\tau, \sigma, \varepsilon \in (0,1]$, such that
\begin{equation}
	\bigl| \bar \phi_h (t)\bigr|_{H^2(\Omega)}
	+
	\bigl| \bar \psi_h (t) \bigr|_{H^2(\Gamma)}
	+
	\bigl| \partial _{\boldsymbol{\nu}} \bar \phi _h(t) \bigr|_{H_\Gamma}
	\le 
	\Lambda_6(t)
	\label{L6}
\end{equation} 
for a.a.\ $t \in (0,T)$.}\smallskip

\paragraph{\bf Proof.} We rewrite (\ref{pLW2}), (\ref{pLW4}), (\ref{pLW6}) in the following elliptic \pier{problem} for $\bar \phi_h$ and $\bar \psi_h$
\begin{equation}
	\begin{cases} 
	 -\Delta \bar\phi_h = \bar\mu_h - \tau \partial _t \hat \phi_h
	- \beta_\varepsilon(\bar \phi_h) - \pi (\bar \phi_h) + \underline{f}\! {_h}
	\quad \mbox{a.e.\ in }\Omega, 
	\\
	(\bar \phi_h)_{|_{\Gamma}} =\bar \psi_h \quad \mbox{a.e.\ on }\Gamma,
	\\ 
	- \Delta_\Gamma \bar \psi_h + \bar \psi_h + \partial _{\boldsymbol{\nu}} \bar \phi_h 
	= \bar w_h - \sigma \partial _t \hat \psi_h 
	- \beta_{\Gamma, \varepsilon}(\bar \psi_h) - \pi_\Gamma(\bar \psi_h) + \underline{g} {_h}
	+\bar \psi_h
	\quad \mbox{a.e.\ on }\Gamma,
	\end{cases}
	\label{el}
\end{equation}
a.e.\ in $(0,T)$. Then it follows from \cite[Lemma~A.1]{MZ05} that the following estimate holds
\begin{align*}
	\bigl| \bar \phi_h \bigr|_{H^2(\Omega)}
	+
	\bigl| \bar \psi_h \bigr|_{H^2(\Gamma)}
	& \le
	C_{\rm MZ} \left( 
	\bigl| \bar\mu_h - \tau \partial _t \hat \phi_h
	- \beta_\varepsilon(\bar \phi_h) - \pi (\bar \phi_h) + \underline{f}\! {_h} \bigr|_H \right. 
	\nonumber \\
	& \quad {} \left. + 
	\bigl| \bar w_h - \sigma \partial _t \hat \psi_h 
	- \beta_{\Gamma, \varepsilon}(\bar \psi_h) - \pi_\Gamma(\bar \psi_h) + \underline{g} {_h}
	+\bar \psi_h
	\bigr|_{H_\Gamma}
	\right) 
\end{align*} 
a.e.\ in $(0,T)$, 
where $C_{\rm MZ}$ is a positive constant independent of 
$h\in(0,h^{**}]$, $\tau, \sigma, \varepsilon \in (0,1]$. 
Accounting for (\ref{M1}), (\ref{L4}) and (\ref{L5}), we see that  
there exist a function $\tilde{\Lambda}_6 \in L^2(0,T)$ and 
a positive constant $\tilde{M}_6$, independent of 
$h\in(0,h^{**}]$, $\tau, \sigma, \varepsilon \in (0,1]$, such that
\begin{equation}
	\bigl| \bar \phi_h (t)\bigr|_{H^2(\Omega)}
	+
	\bigl| \bar \psi_h (t) \bigr|_{H^2(\Gamma)}
	\le 
	\tilde{M}_6 
	\left( \tilde{\Lambda}_6(t)+\bigl| \partial _{\boldsymbol{\nu}} \bar \phi _h(t) \bigr|_{H_\Gamma} \right)
	\label{tL6}
\end{equation}
for a.a.\ $t \in (0,T)$. We observe that, for some $3/2<s<2$, there exists a positive constant $C_7$ such that
\begin{equation}
	\bigl| \partial_{\boldsymbol{\nu}} \bar \phi_h(t)\bigr|_{H_\Gamma} 
	\le C_7 \bigl| \bar \phi_h (t) \bigr|_{H^s(\Omega)} 
	\label{trace}
\end{equation}
for a.a.\ $t \in (0,T)$ (see, e.g., \cite[Theorem~2.25, p.~1.62]{BG87}), and consequently, by the compactness inequality 
(see, e.g., \cite[Lemme~5.1, p.~58]{Lio69} or \cite[Theorem~16.4, p.~102]{LM70}) as well as 
(\ref{tL6}) and (\ref{trace}), we infer 
\begin{align*}
	&\bigl| \bar \phi_h (t)\bigr|_{H^2(\Omega)}
	+
	\bigl| \bar \psi_h (t)\bigr|_{H^2(\Gamma)}
	+
	\bigl| \partial _{\boldsymbol{\nu}} \bar \phi _h(t) \bigr|_{H_\Gamma}
	\nonumber\\
	& \quad \le
	\tilde{M}_6 
	\left( \tilde{\Lambda}_6(t)+C_7 \bigl| \bar \phi_h (t) \bigr|_{H^s(\Omega)}  \right)
	+ C_7 \bigl| \bar \phi_h (t) \bigr|_{H^s(\Omega)} 
	\nonumber \\
	& \quad \le \frac{1}{2}\bigl| \bar \phi_h (t)\bigr|_{H^2(\Omega)}+ 
	C_8 \left( \tilde{\Lambda}_6(t)+ \bigl| \bar \phi_h (t)\bigr|_{V} \right)
\end{align*} 
for a.a.\ $t \in (0,T)$, whence, also by (\ref{M1}), we conclude \pier{that (\ref{L6}) holds for}
some function $\Lambda_6 \in L^2(0,T)$, independent of 
$h\in(0,h^{**}]$, $\tau, \sigma, \varepsilon \in (0,1]$.  
\hfill $\Box$\medskip

Now, using these uniform estimates, we can discuss the existence of weak solutions to the viscous problem for the original system \eqref{LW1}--\eqref{LW7}, by taking $h\to 0$ and $\varepsilon \to 0$. 
The subscripts of $\tau$ and $\sigma$ for functions mean 
the dependence on parameters $\tau, \sigma \in (0,1]$, however 
in the next proposition we omit them for simplicity.\medskip

\paragraph{\bf Proposition~3.2}
{\it Assume the (A1)--(A5) are satisfied. For each $\tau, \sigma \in (0,1]$, there exists a 
sextuplet $(\phi, \mu, \xi, \psi, w, \zeta)$ of functions  
\begin{align*}
	& \phi :=\phi_{\tau,\sigma} \in H^1(0,T;V^*)\cap L^\infty (0,T;V) \cap L^2 \bigl(0,T;H^2(\Omega) \bigr), \\
	& \partial_{\boldsymbol{\nu}} \phi \in L^2(0,T;H_\Gamma), \quad 
	\tau \phi \in H^1(0,T;H), \\
	& \mu:=\mu_{\tau,\sigma} \in L^2(0,T;V), \quad \tau \mu \in L^2(0,T;W), \quad 
	\xi :=\xi_{\tau,\sigma} \in L^2(0,T;H), \\
	& \psi :=\psi _{\tau,\sigma} \in H^1(0,T;V_\Gamma^*)\cap L^\infty (0,T;V_\Gamma) \cap L^2 \bigl( 0,T;H^2(\Gamma) \bigr), 
	\quad \sigma \psi \in H^1(0,T;H_\Gamma),\\
	& w:=w_{\tau,\sigma} \in L^2(0,T;V_\Gamma), \quad \sigma w \in L^2 \bigl( 0,T;H^2(\Gamma) \bigr), \quad 
	\zeta:=\zeta_{\tau,\sigma} \in L^2(0,T;H_\Gamma)
\end{align*} 
such that they satisfy the following viscous {C}ahn--{H}illiard system:
\begin{align} 
	& \partial_t \phi - \Delta \mu =0 
	& \mbox{a.e.\ in }Q, 
	\label{visLW1}
	\\
	& \mu = \tau \partial_t \phi -\Delta \phi + \xi + \pi (\phi) -f,\quad \xi \in \beta (\phi)
	& \mbox{a.e.\ in }Q, 
	\label{visLW2}
	\\
	& \partial_{\boldsymbol{\nu}} \mu =0 
	& \mbox{a.e.\ on }\Sigma, 
	\label{visLW3}
	\\
	& \phi_{|_{\Gamma}} =\psi 
	& \mbox{a.e.\ on }\Sigma,
	\label{visLW4}
	\\ 
	& \partial_t \psi -\Delta_\Gamma w = 0 
	& \mbox{a.e.\ on }\Sigma, 
	\label{visLW5}
	\\
	& w = \sigma \partial_t \psi + \partial _{\boldsymbol{\nu}} \phi - \Delta_\Gamma \psi 
	+ \zeta + \pi_\Gamma(\psi)-g,\quad 
	\zeta \in \beta_\Gamma (\psi)
	& \mbox{a.e.\ on }\Sigma,
	\label{visLW6}
	\\
	&\phi(0)=\phi_0
	\quad 
	\mbox{a.e.\ in }\Omega, 
	\quad 
	\psi (0)=\psi_0
	\quad \mbox{a.e.\ on }\Gamma. &
	\label{visLW7}
\end{align}
}

\paragraph{\bf Proof.} 
Let $\tau, \sigma \in (0,1]$. 
Recalling (\ref{M1}), (\ref{L4}),  (\ref{L5}) and (\ref{L6}), we deduce that there 
exist positive constants $M_6, M_7$, independent of 
$h\in(0,h^{**}]$, $\tau, \sigma, \varepsilon \in (0,1]$, such that
\begin{align}
	& | \bar{\mu}_h |_{L^2(0,T;V)} 
	+ 
	| \bar{w}_h |_{L^2(0,T;V_\Gamma)} 
	+
	\bigl| \beta_\varepsilon \bigl( \bar \phi_h \bigr) \bigr|_{L^2(0,T;H)}
	+
	\bigl| \beta_{\Gamma, \varepsilon} \bigl( \bar \psi_h \bigr) \bigr|_{L^2(0,T;H_\Gamma)}
	\nonumber \\
	& \quad 
	+ \bigl| \bar \phi_h \bigr|_{L^2(0,T;H^2(\Omega))}
	+
	\bigl| \bar \psi_h \bigr|_{L^2(0,T;H^2(\Gamma))}
	+
	\bigl| \partial _{\boldsymbol{\nu}} \bar \phi _h \bigr|_{L^2(0,T;H_\Gamma)}
	\le M_6,
	\label{M6}
	\\[0.2cm]
	& \sqrt{\tau} | \bar \mu _h |_{L^2(0,T;W)} + \sqrt{\sigma}| \bar w _h |_{L^2 (0,T;H^2(\Gamma))} \le M_7,
	\label{M7}
\end{align}
for all $h \in(0,h^{**}]$ and $\varepsilon \in (0,1]$\pier{, where also a comparison in 
(\ref{pLW1}), (\ref{pLW3}) and (\ref{pLW5}) has been used for \eqref{M7}.}
Next, using (\ref{tool1})--(\ref{tool3}), we observe that 
\begin{align*}
	& \bigl| \hat \phi _h \bigr|_{L^\infty(0,T;V)} \le 
	\max\bigl\{ |\phi_0|_V, \bigl|\bar\phi_h \bigr|_{L^\infty(0,T;V)} \bigr\}\le \sqrt{M_1}+|\phi_0|_V, 
	\\
	& \tau \bigl| \hat\phi_h - \bar\phi_h\bigr|_{L^2(0,T;H)}^2
	= \frac{h^2}{3} \tau \bigl|\partial _t \hat \phi_h \bigr|_{L^2(0,T;H)}^2 \le \frac{h^2}{3}M_1 \le \frac{1}{3}M_1, 
	\\
	& \bigl| h \partial_t \hat \mu_h \bigr|_{L^2(0,T;H)} 
	+ \bigl| \partial _t \hat \phi_h \bigr|_{L^2(0,T;V^*)} \le 3 \sqrt{M_1}, 
	\\
	& \bigl|\hat \mu_h \bigr|_{L^2(0,T;V)}^2 
	\le \frac{h}{2} |\mu_0|_V^2+|\bar \mu_h |_{L^2(0,T;V)}^2 
	\le M_6^2,\\
	& \bigl| \hat \psi _h \bigr|_{L^\infty(0,T;V_\Gamma)} \le 
	\max\bigl\{ |\psi_0|_{V_\Gamma}, \bigl|\bar\psi_h \bigr|_{L^\infty(0,T;V_\Gamma)} \bigr\}
	\le \sqrt{M_1}+|\psi_0|_{V_\Gamma}, 
	\\
	& \sigma \bigl| \hat\psi_h - \bar\psi_h \bigr|_{L^2(0,T;H_\Gamma)}^2
	= \frac{h^2}{3} \sigma \bigl|\partial _t \hat \psi_h \bigr|_{L^2(0,T;H_\Gamma)}^2 
	\le \frac{h^2}{3}M_1 	
	\le \frac{1}{3}M_1, 
	\\
	& \bigl| h \partial_t \hat w_h \bigr|_{L^2(0,T;H_\Gamma)} 
	+ \bigl| \partial _t \hat \psi_h \bigr|_{L^2(0,T;V_\Gamma^*)} \le 3 \sqrt{M_1}, 
	\\
	& \bigl|\hat w_h \bigr|_{L^2(0,T;V_\Gamma)}^2 
	\le \frac{h}{2} |w_0|_{V_\Gamma}^2+|\bar w_h |_{L^2(0,T;V_\Gamma)}^2 
	\le M_6^2,
\end{align*}
for all $h \in(0,h^{**}]$ and $\varepsilon \in (0,1]$.  
Then, there exist functions $\phi_\varepsilon, \mu_\varepsilon, \psi_\varepsilon, w_\varepsilon$ and 
a subsequence $\{h_k \}_{k \in \mathbb{N}}$ of $h \to 0$ such that 
\begin{align*}
	\hat \phi_{h_k} \to \phi_\varepsilon 
	& \quad \mbox{weakly star in } H^1(0,T;V^*) \cap L^\infty (0,T;V)\\
	& \quad \mbox{and strongly in } C\bigl( [0,T]; H \bigr), \\
	\pier{\tau \hat \phi_{h_k} \to \tau \phi_\varepsilon}
	& \quad \pier{\mbox{weakly in } H^1(0,T;H),} \\	
	\bar \phi_{h_k} \to \phi_\varepsilon 
	& \quad \mbox{weakly star in } L^\infty(0,T;V) \cap L^2 \bigl(0,T;H^2(\Omega) \bigr) \\
	& \quad \mbox{and strongly in } L^2(0,T; H), \\
	\bar \mu_{h_k} \to \mu _\varepsilon
	& \quad \mbox{weakly in } L^2 (0,T;W), \\
	h_k \hat \mu_{h_k} \to 0
	& \quad \mbox{weakly in } H^1 (0,T;H) \ \mbox{and strongly in } L^2(0,T;V), \\
	\hat \psi_{h_k} \to \psi_\varepsilon 
	& \quad \mbox{weakly star in } H^1(0,T;V_\Gamma^*) \cap L^\infty (0,T;V_\Gamma) \\
	& \quad \mbox{and strongly in } C\bigl( [0,T]; H_\Gamma \bigr), \\
		\pier{\sigma \hat \psi_{h_k} \to \sigma \psi_\varepsilon} 
	& \quad \pier{\mbox{weakly in } H^1(0,T;H_\Gamma),} \\
	\bar \psi_{h_k} \to \psi_\varepsilon 
	& \quad \mbox{weakly star in } L^\infty(0,T;V_\Gamma) \cap L^2 \bigl(0,T;H^2(\Gamma) \bigr) \\
	& \quad \mbox{and strongly in } L^2(0,T; H_\Gamma), \\
	\bar w_{h_k} \to w_\varepsilon
	& \quad \mbox{weakly in } L^2 \bigl(0,T;H^2(\Gamma) \bigr), \\
	h_k \hat w_{h_k} \to 0
	& \quad \mbox{weakly in } H^1 (0,T;H_\Gamma) \ \mbox{and strongly in } L^2(0,T;V_\Gamma), \\
	\partial _{\boldsymbol{\nu}} \bar \phi_{h_k} \to \partial _{\boldsymbol{\nu}} \phi_\varepsilon 
	& \quad \mbox{weakly in } L^2 (0,T;H_\Gamma)
\end{align*}
as $k\to +\infty$, where we applied the compactness results \cite[Section 8, Corollary 4]{Sim87} 
\pier{and \eqref{tool3} to obtain the} 
strong convergences. \pier{We recall that $\tau$ and $\sigma$ are positive and fixed in this limit procedure.}
Moreover, due to the {L}ipschitz continuity of 
$\beta_\varepsilon$, $\pi$, $\beta_{\Gamma, \varepsilon}$ 
and $\pi_\Gamma$, we have that
\begin{align*}
	 \beta_\varepsilon \bigl( \bar \phi_{h_k} \bigr) \to \beta_\varepsilon(\phi_\varepsilon), 
	\quad 
	\pi \bigl( \bar \phi_{h_k} \bigr) \to \pi(\phi_\varepsilon) 
	& \quad \mbox{strongly in } L^2(0,T; H), \\
	 \beta_{\Gamma,\varepsilon} \bigl( \bar \psi_{h_k} \bigr) \to \beta_{\Gamma,\varepsilon}(\psi_\varepsilon), 
	\quad 
	\pi_\Gamma \bigl( \bar \psi_{h_k} \bigr) \to \pi_\Gamma (\psi_\varepsilon) 
	& \quad \mbox{strongly in } L^2(0,T; H_\Gamma)
\end{align*}
as $k\to +\infty$. 
Besides, \pier{it is not difficult to check that
\begin{align*}
	f_h \to f
	& \quad \mbox{strongly in } L^2(0,T; V), \\
	g_h \to g
	& \quad \mbox{strongly in } L^2(0,T; V_\Gamma)
\end{align*}
as $h \to 0$  (in \cite[Appendix]{CK19} the argument is fully detailed)}. 
Based on all these convergence results, we can pass to the limit as $h_k \to 0$ in problem 
(\ref{pLW1})--(\ref{pLW7}) and find that the 
quadruplet $(\phi_\varepsilon, \mu_\varepsilon, \psi_\varepsilon, w_\varepsilon)$ solves 
\begin{align} 
	& \partial_t \phi_\varepsilon - \Delta \mu_\varepsilon =0 
	& \mbox{a.e.\ in }Q, 
	\label{yLW1}
	\\
	& \mu_\varepsilon = \tau \partial_t \phi_\varepsilon -\Delta \phi_\varepsilon + \beta_\varepsilon(\phi_\varepsilon) + \pi (\phi_\varepsilon) -f
	& \mbox{a.e.\ in }Q, 
	\label{yLW2}
	\\
	& \partial_{\boldsymbol{\nu}} \mu_\varepsilon =0 
	& \mbox{a.e.\ on }\Sigma, 
	\label{yLW3}
	\\
	& (\phi_\varepsilon )_{|_{{\Gamma}}} =\psi_\varepsilon 
	& \mbox{a.e.\ on }\Sigma,
	\label{yLW4}
	\\ 
	& \partial_t \psi_\varepsilon -\Delta_\Gamma w_\varepsilon = 0 
	& \mbox{a.e.\ on }\Sigma, 
	\label{yLW5}
	\\
	& w_\varepsilon=\partial_{\boldsymbol{\nu}} \phi_\varepsilon+\sigma \partial _t \psi _\varepsilon
	-\Delta_\Gamma 
	\psi_\varepsilon+\beta_{\Gamma,\varepsilon}(\psi_\varepsilon)+\pi_{\Gamma} (\psi_\varepsilon)-g
	& \mbox{a.e.\ on }\Sigma, 
	\label{yLW6}
	\\
	& \phi_\varepsilon(0)=\phi_0
	\quad 
	\mbox{a.e.\ in }\Omega, 
	\quad 
	\psi _\varepsilon(0)=\psi_0
	\quad \mbox{a.e.\ on }\Gamma. &
	\label{yLW7}
\end{align}
Moreover, by weak or weak star lower semicontinuity of norms, we see that the following estimates 
hold (see (\ref{M1}), (\ref{M6}), (\ref{M7}))
\begin{align}
	& | \partial_t \phi_\varepsilon |_{L^2(0,T;V^*)} 
	\le \liminf_{k \to \infty} \bigl| \partial _t \hat \phi_{h_k} \bigr|_{L^2(0,T;V^*)} 
	\le 3\sqrt{M_1}, 
	\label{est1}\\
	& | \phi_\varepsilon |_{L^\infty(0,T;V)} 
	\le \liminf_{k \to \infty} \bigl| \hat \phi_{h_k} \bigr|_{L^\infty (0,T;V)} 
	\le \sqrt{M_1}+|\phi_0|_V, 
	\label{est2}\\
	& | \phi_\varepsilon |_{L^2(0,T;H^2(\Omega))} 
	\le \liminf_{k \to \infty} \bigl| \bar \phi_{h_k} \bigr|_{L^2 (0,T;H^2(\Omega))}
	\le M_6, 
	\label{est3}\\
	& \sqrt{\tau} | \partial_t \phi_\varepsilon |_{L^2(0,T;H)} 
	\le \liminf_{k \to \infty} \sqrt{\tau} | \partial_t \hat \phi_{h_k} |_{L^2(0,T;H)} 
	\le \sqrt{M_1}, 
	\label{est4}\\
	& |\mu_\varepsilon|_{L^2(0,T;V)} \le 
	 \liminf_{k \to \infty} | \bar \mu_{h_k} |_{L^2(0,T;V)} \le M_6, 
	 \label{est5}\\
	 & \sqrt{\tau} |\mu_\varepsilon|_{L^2(0,T;W)} \le 
	 \liminf_{k \to \infty}  \sqrt{\tau} | \bar \mu_{h_k} |_{L^2(0,T;W)} \le M_7, 
	 \label{est6}\\
	& | \partial_t \psi_\varepsilon |_{L^2(0,T;V_\Gamma^*)} 
	\le \liminf_{k \to \infty} \bigl| \partial_t \hat \psi_{h_k} \bigr|_{L^2(0,T;V_\Gamma^*)} 
	\le 3\sqrt{M_1}, 
	\label{est7}\\
	& | \psi_\varepsilon |_{L^\infty(0,T;V_\Gamma)} 
	\le \liminf_{k \to \infty} \bigl| \hat \psi_{h_k} \bigr|_{L^\infty (0,T;V_\Gamma)} 
	\le \sqrt{M_1}+|\psi_0|_{V_\Gamma}, 
	\label{est8}\\
	& | \psi_\varepsilon |_{L^2(0,T;H^2(\Gamma))} 
	\le \liminf_{k \to \infty} \bigl| \bar \psi_{h_k} \bigr|_{L^2 (0,T;H^2(\Gamma))}
	\le M_6, 
	\label{est9}\\
	& \sqrt{\sigma} | \partial_t \psi_\varepsilon |_{L^2(0,T;H_\Gamma)} 
	\le \liminf_{k \to \infty} \sqrt{\sigma} \bigl| \partial_t \hat \psi_{h_k} \bigr|_{L^2(0,T;H_\Gamma)}
	 \le \sqrt{M_1}, 
	 \label{est10}\\
	& |w_\varepsilon|_{L^2(0,T;V_\Gamma)} \le 
	 \liminf_{k \to \infty} | \bar w_{h_k} |_{L^2(0,T;V_\Gamma)} \le M_6, 
	 \label{est11}\\
	 & \sqrt{\sigma} |w_\varepsilon|_{L^2(0,T;H^2(\Gamma))} \le 
	 \liminf_{k \to \infty}  \sqrt{\sigma} | \bar w_{h_k} |_{L^2(0,T;H^2(\Gamma))} \le M_7, 
	 \label{est12}\\
	 & |\partial_{\boldsymbol{\nu}} \phi_\varepsilon|_{L^2(0,T;H_\Gamma)}
	 \le  \liminf_{k \to \infty} \bigl|\partial_{\boldsymbol{\nu}} \bar \phi_{h_k} \bigr|_{L^2(0,T;H_\Gamma)} \le M_6, 
	 \label{est13}\\
	 & \bigl| \beta_\varepsilon ( \phi_\varepsilon ) \bigr|_{L^2(0,T;H)}
	 = \lim_{k \to \infty}	 \bigl| \beta_\varepsilon \bigl( \bar\phi_{h_k} \bigr) \bigr|_{L^2(0,T;H)} \le M_6, 
	 \label{est14}\\
	 & \bigl| \beta_{\Gamma,\varepsilon} ( \psi_\varepsilon) \bigr|_{L^2(0,T;H_\Gamma)}
	 = \lim_{k \to \infty}	 \bigl| \beta_{\Gamma, \varepsilon} \bigl( \bar\psi_{h_k} \bigr) \bigr|_{L^2(0,T;H_\Gamma)} \le M_6.
	 \label{est15}
\end{align}
Due to the uniform estimates (\ref{est1})--(\ref{est15}), 
we are able to pass to the limit along a subsequence $\{\varepsilon_k\}_{k \in\mathbb{N}}$ of $\varepsilon$, 
in the problem (\ref{yLW1})--(\ref{yLW7}) by finding elements 
$$
\phi:=\phi_{\tau,\sigma},\ \ \mu:=\mu_{\tau,\sigma},\ \ \xi:=\xi_{\tau,\sigma}, \ \ 
\psi:=\psi_{\tau,\sigma},\ \ w:=w_{\tau,\sigma},\ \ \zeta:=\zeta_{\tau,\sigma}
$$ 
such that 
\begin{align*}
	\phi_{\varepsilon_k} \to \phi 
	& \quad \mbox{weakly star in } H^1(0,T;V^*) \cap L^\infty (0,T;V) \cap L^2 \bigl(0,T;H^2(\Omega) \bigr)\\
	& \quad \mbox{and strongly in } C\bigl( [0,T]; H \bigr) \cap L^2(0,T;V), \\
	\tau \phi_{\varepsilon_k} \to \tau \phi 
	& \quad \mbox{weakly in } H^1(0,T;H), \\
	\mu_{\varepsilon_k} \to \mu
	& \quad \mbox{weakly in } L^2 (0,T;W), \\
	\beta_{\varepsilon_k} ( \phi_{\varepsilon_k} ) \to \xi 
	& \quad \mbox{weakly in } L^2(0,T; H), \\
	\psi_{\varepsilon_k} \to \psi 
	& \quad \mbox{weakly star in } H^1(0,T;V_\Gamma^*) \cap L^\infty (0,T;V_\Gamma) \cap L^2 \bigl(0,T;H^2(\Gamma) \bigr)\\
	& \quad \mbox{and strongly in } C\bigl( [0,T]; H_\Gamma \bigr) \cap L^2(0,T;V_\Gamma), \\
	\sigma \psi_{\varepsilon_k} \to \sigma \psi 
	& \quad \mbox{weakly in } H^1(0,T;H_\Gamma), \\
	w_{\varepsilon_k} \to w
	& \quad \mbox{weakly in } L^2 \bigl(0,T;H^2(\Gamma) \bigr), \\
	\beta_{\Gamma,\varepsilon_k} ( \psi_{\varepsilon_k} ) \to \zeta
	& \quad \mbox{weakly in } L^2(0,T; H_\Gamma), \\
	\partial _{\boldsymbol{\nu}} \phi_{\varepsilon_k} \to \partial _{\boldsymbol{\nu}} \phi 
	& \quad \mbox{weakly in } L^2 (0,T;H_\Gamma)
\end{align*}
as $k\to +\infty$, due to the compactness theorems again.
Now, we 
observe that 
$\xi \in \beta (\phi)$ a.e.\ in $Q$ and 
$\zeta \in \beta_\Gamma (\psi)$ a.e.\ on $\Sigma$, 
due to the maximal monotonicity of $\beta$ and $\beta_\Gamma$, 
and the weak-strong convergence for $\beta_{\varepsilon_k}(\phi_{\varepsilon_k})$ and 
$\phi_{\varepsilon_k}$ in $L^2(0,T;H)$, and for 
$\beta_{\Gamma,\varepsilon_k}(\psi_{\varepsilon_k})$ and 
$\psi_{\varepsilon_k}$ in $L^2(0,T;H_\Gamma)$, respectively. 
Finally, we observe that 
\begin{align*}
	\pi(\phi_{\varepsilon_k}) \to \pi (\phi) 
	& \quad 
	\mbox{strongly in }  C\bigl( [0,T]; H \bigr), \\
	\pi_\Gamma (\psi_{\varepsilon_k}) \to \pi_\Gamma (\psi) 
	& \quad 
	\mbox{strongly in }  C\bigl( [0,T]; H_\Gamma \bigr). 
\end{align*} 
Thus, we can pass to the limit \pier{as} $k \to \infty$ in the regularized problem (\ref{yLW1})--(\ref{yLW7}) to obtain the viscous Cahn--Hilliard system 
(\ref{visLW1})--(\ref{visLW7}). \hfill $\Box$\smallskip

\subsection{From viscous to pure {C}ahn--{H}illiard system}

As a summary of the previous subsection, we can find a 
sextuplet $(\phi_{\tau,\sigma}, \mu_{\tau,\sigma}, \xi_{\tau,\sigma}, \psi_{\tau,\sigma}, w_{\tau,\sigma}, \zeta_{\tau,\sigma})$ of functions, \pier{depending} on $\tau, \sigma \in (0,1]$, 
such that it satisfies the viscous {C}ahn--{H}illiard system (\ref{visLW1})--(\ref{visLW7}). 
Moreover, in (\ref{est1})--(\ref{est15}), from weak or weak star lower semicontinuity of norms, we also know that 
there exists a positive constant $M_8$, independent of  $\tau, \sigma \in (0,1]$, such that 
\begin{align}
	& | \phi_{\tau,\sigma} |_{H^1(0,T;V^*)} 
	+
	| \phi_{\tau,\sigma} |_{L^\infty(0,T;V)} 
	+
	| \phi_{\tau,\sigma} |_{L^2(0,T;H^2(\Omega))} 
	+
	\sqrt{\tau} | \phi_{\tau,\sigma} |_{H^1(0,T;H)} 
	\nonumber \\
	&\quad  +
	|\mu_{\tau,\sigma}|_{L^2(0,T;V)} 
	+
	 \sqrt{\tau} |\mu_{\tau,\sigma}|_{L^2(0,T;W)} 
	+
	| \psi_{\tau,\sigma} |_{H^1(0,T;V_\Gamma^*)} 
	+
	| \psi_{\tau,\sigma} |_{L^\infty(0,T;V_\Gamma)} 
	\nonumber  \\
	&\quad 
	+
	| \psi_{\tau,\sigma} |_{L^2(0,T;H^2(\Gamma))} 
	+
	\sqrt{\sigma} | \psi_{\tau,\sigma} |_{H^1(0,T;H_\Gamma)} 
	+
	|w_{\tau,\sigma} |_{L^2(0,T;V_\Gamma)} 
	+
	 \sqrt{\sigma} |w_{\tau,\sigma} |_{L^2(0,T;H^2(\Gamma))}
	 \nonumber \\
	&\quad  
	+
	 |\partial_{\boldsymbol{\nu}} \phi_{\tau,\sigma}|_{L^2(0,T;H_\Gamma)}
	+
	 | \xi_{\tau,\sigma} |_{L^2(0,T;H)}
	+
	 | \zeta _{\tau,\sigma} |_{L^2(0,T;H_\Gamma)}
	 \le M_8.
\end{align} 

Now we are in a position to prove Theorem~2.1.\smallskip

\paragraph{\bf Proof of Theorem~2.1.}
We obtain from (\ref{visLW1}), (\ref{visLW3}) and (\ref{visLW5}), 
the following variational \pier{formulations}:
\begin{align}
	& \langle \partial _t \phi_{\tau,\sigma}, z \rangle_{V^*,V} + 
	\int_\Omega \nabla \mu_{\tau,\sigma} \cdot \nabla z dx =0 
	\quad \mbox{for all } z \in V,  
	\label{weakf1}
	\\
	& \langle \partial _t \psi_{\tau,\sigma}, z_\Gamma \rangle_{V^*,V} + 
	\int_\Gamma \nabla_\Gamma w_{\tau,\sigma} \cdot \nabla_\Gamma z_\Gamma d\Gamma =0 
	\quad \mbox{for all } z_\Gamma \in V_\Gamma,
	\label{weakf2}
\end{align}
a.e.\ in $(0,T)$. At this point, we can pass to the limit as 
either $\tau \to 0$ or $\sigma \to 0$, or both $\tau, \sigma \to 0$
in order to obtain a partially viscous {C}ahn--{H}illiard system or a 
pure {C}ahn--{H}illiard system at the limit. 
Let us detail only the last case with $(\tau, \sigma) \to (0,0)$ 
along a joint subsequence $(\tau_k, \sigma_k)$. We see that 
there exists a sextuplet $(\phi, \mu, \xi, \psi, w, \zeta)$ such that 
\begin{align*}
	\phi_{\tau_k,\sigma_k} \to \phi
	& \quad \mbox{weakly star in } H^1(0,T;V^*) \cap L^\infty (0,T;V) \cap L^2 \bigl(0,T;H^2(\Omega) \bigr)\\
	& \quad \mbox{and strongly in } C\bigl( [0,T]; H \bigr) \cap L^2(0,T;V), \\
	\tau_k \phi_{\tau_k,\sigma_k} \to 0
	& \quad \mbox{strongly in } H^1(0,T;H), \\
	\mu_{\tau_k, \sigma_k} \to \mu
	& \quad \mbox{weakly in } L^2 (0,T;V), \\
	\xi_{\tau_k, \sigma_k} \to \xi 
	& \quad \mbox{weakly in } L^2(0,T; H), \\
	\psi_{\tau_k,\sigma_k} \to \psi 
	& \quad \mbox{weakly star in } H^1(0,T;V_\Gamma^*) \cap L^\infty (0,T;V_\Gamma) \cap L^2 \bigl(0,T;H^2(\Gamma) \bigr)\\
	& \quad \mbox{and strongly in } C\bigl( [0,T]; H_\Gamma \bigr) \cap L^2(0,T;V_\Gamma), \\
	\sigma_k \psi_{\tau_k,\sigma_k} \to 0 
	& \quad \mbox{strongly in } H^1(0,T;H_\Gamma), \\
	w_{\tau_k,\sigma_k} \to w
	& \quad \mbox{weakly in } L^2 (0,T;V_\Gamma ), \\
	\zeta_{\tau_k,\sigma_k} \to \zeta
	& \quad \mbox{weakly in } L^2(0,T; H_\Gamma), \\
	\partial _{\boldsymbol{\nu}} \phi_{\tau_k,\sigma_k} \to \partial _{\boldsymbol{\nu}} \phi 
	& \quad \mbox{weakly in } L^2 (0,T;H_\Gamma)
\end{align*}
as $k \to +\infty$. 
Different from the previous subsection, we can pass to the limit in
\begin{equation*}
	\xi_{\tau_k, \sigma_k} \in \beta (\phi_{\tau_k, \sigma_k}) \quad \mbox{a.e.\ in } Q, \quad 
	\zeta_{\tau_k, \sigma_k} \in \beta_\Gamma (\psi_{\tau_k, \sigma_k}) \quad \mbox{a.e.\ in } \Sigma
\end{equation*}
just using the demi-closedness of $\beta$ and $\beta_\Gamma$, respectively, to obtain the same inclusions at the limit. 
To complete this limiting procedure, we pass to the limit in (\ref{visLW2}), (\ref{visLW4}), (\ref{visLW6}), (\ref{visLW7}), 
(\ref{weakf1}), (\ref{weakf2}) 
to obtain (\ref{weak1})--(\ref{weak6})\pier{. Hence, we arrive at} the conclusion. \hfill $\Box$\medskip

Let us remark that, if we let only one of the parameters $\tau, \sigma$ go to $0$, \pier{then} 
we also have the convergence 
\begin{equation*}
	\mu_{\tau, \sigma_k} \to \mu_\tau 
	\quad \mbox{weakly in } L^2 (0,T;W) 
	\quad \mbox{if only } \sigma_k \to 0, 
\end{equation*}
or the convergence 
\begin{equation*}
	w_{\tau_k,\sigma} \to w_\sigma
	\quad \mbox{weakly in } L^2 \bigl(0,T;H^2(\Gamma) \bigr)
	\quad \mbox{if only } \tau_k \to 0.
\end{equation*}
In these cases, we can keep the smoothness of the time derivative, more precisely, 
$\phi_\tau \in H^1(0,T;H)$ if $\tau>0$, or 
$\psi _\sigma \in H^1(0,T;H_\Gamma)$ if $\sigma>0$, respectively.\medskip

\paragraph{\bf Proof of Theorem~2.2.}
We now prove a continuous dependence estimates with two weak solutions 
\begin{equation*}
		\bigl( \phi^{(i)}, \mu^{(i)}, \xi^{(i)}, \psi^{(i)}, w^{(i)}, \zeta^{(i)} \bigr) \quad 
		\mbox{for } i=1,2,
\end{equation*}
corresponding to the initial data 
\begin{equation*}
	\bigl( \phi_0^{(i)}, \psi_0^{(i)} \bigr) \in V^* \times V_\Gamma^* \quad 
	\mbox{for } i=1,2,
\end{equation*}
\pier{satisfying (\ref{pier3}), (\ref{pier4}),}  and the sources 
\begin{equation*}
	\bigl( f^{(i)}, g^{(i)} \bigr) \in L^2(0,T;V^*) \times L^2(0,T;V_\Gamma^*) \quad 
	\mbox{for } i=1,2.
\end{equation*}
We take the difference of (\ref{weak1}) and choose $z:=1$ to obtain that 
\begin{equation*}
	\bigl\langle \partial_t \bigl(\phi^{(1)}-\phi^{(2)}\bigr) ,1 \bigr\rangle_{V^*,V}=0 
\end{equation*}
a.e.\ in $(0,T)$, whence
\begin{equation*}
	\bigl\langle \phi^{(1)}(t)-\phi^{(2)}(t), 1 \bigr\rangle_{V^*,V}
	= 
	\bigl\langle \phi^{(1)}_0-\phi^{(2)}_0, 1 \bigr\rangle_{V^*,V}
	= 0
\end{equation*}
for all $t \in [0,T]$. 
Then, we can take $z:={\mathcal N}_\Omega (\phi^{(1)}-\phi^{(2)})$ as a 
test function in the difference of (\ref{weak1}) and obtain
\begin{equation}
	\frac{1}{2}
	\bigl| \phi^{(1)}(t)-\phi^{(2)}(t)\bigr|_{V_{0*}}^2
	+ \int_0^t\!\!\int_\Omega 
	\bigl( \mu^{(1)}-\mu^{(2)} \bigr) 
	\bigl( \phi^{(1)}-\phi^{(2)} \bigr) dxds 
	=
	\frac{1}{2}
	\bigl| \phi^{(1)}_0-\phi^{(2)}_0 \bigr|_{V_{0*}}^2
	\label{conti1}
\end{equation}
for all $t \in [0,T]$. 
By operating on the difference of (\ref{weak4}) in the same way, that is, $z_\Gamma:=1$ first and 
$z_\Gamma:={\mathcal N}_\Gamma (\psi^{(1)}-\psi^{(2)})$ second, we obtain the similar formula
\begin{equation}
	\frac{1}{2}
	\bigl| \psi^{(1)}(t)-\psi^{(2)}(t)\bigr|_{V_{\Gamma,0*}}^2
	+ \int_0^t\!\!\int_\Gamma
	\bigl( w^{(1)}-w^{(2)} \bigr) 
	\bigl( \psi^{(1)}-\psi^{(2)} \bigr) d\Gamma ds 
	=
	\frac{1}{2}
	\bigl| \psi^{(1)}_0-\psi^{(2)}_0 \bigr|_{V_{\Gamma,0*}}^2
	\label{conti2}
\end{equation}
for all $t \in [0,T]$. 
Next, we multiply the difference of \pier{the equalities in} (\ref{weak2}) by $\phi^{(1)}-\phi^{(2)}$ and 
integrate the resultant with respect to space and time. 
Using (\ref{weak3}) and (\ref{weak5}), \pier{we infer that}
\begin{align*}
	& \int_0^t\!\! \int_\Omega 
	\bigl( \mu^{(1)}-\mu^{(2)} \bigr) 
	\bigl( \phi^{(1)}-\phi^{(2)} \bigr) dxds 
	+
	\int_0^t\!\! \int_\Gamma
	\bigl( w^{(1)}-w^{(2)} \bigr) 
	\bigl( \psi^{(1)}-\psi^{(2)}\bigr) d\Gamma ds 
	\nonumber \\
	& = \int_0^t\!\!\int _\Omega \bigl| \nabla \bigl(\phi^{(1)}-\phi^{(2)} \bigr) \bigr|^2dx ds 
	+ \int_0^t\!\!\int _\Gamma \bigl| \nabla_\Gamma \bigl(\psi^{(1)}-\psi^{(2)} \bigr) \bigr|^2d\Gamma ds 
	\nonumber \\
	& \quad {} + \int_0^t\!\!\int _\Omega \bigl( \xi ^{(1)}-\xi ^{(2)} \bigr) 
	\bigl( \phi^{(1)}-\phi^{(2)} \bigr)dx ds 
	+ \int_0^t\!\!\int _\Gamma \bigl( \zeta^{(1)}-\zeta^{(2)} \bigr) 
	\bigl( \psi^{(1)}-\psi^{(2)}\bigr)  d\Gamma ds 
	\nonumber \\
	& \quad {}+ \int_0^t\!\!\int _\Omega \bigl( \pi\bigl(\phi^{(1)}\bigr)-\pi\bigl(\phi^{(2)}\bigr) \bigr)
	\bigl( \phi^{(1)}-\phi^{(2)} \bigr)dx ds 
	- \int_0^t \bigl\langle 
	f^{(1)}-f^{(2)}, \phi^{(1)}-\phi^{(2)} \bigr\rangle_{V^*,V} ds
	\nonumber \\
	& \quad {}
	+ \int_0^t\!\!\int _\Gamma \bigl( \pi_\Gamma\bigl(\psi^{(1)}\bigr)-\pi_\Gamma\bigl(\psi^{(2)}\bigr) \bigr)
	\bigl( \psi^{(1)}-\psi^{(2)}\bigr)  d\Gamma ds 
	- \int_0^t \bigl\langle g^{(1)}-g^{(2)},
	\psi^{(1)}-\psi^{(2)} \bigr\rangle_{V_\Gamma^*,V_\Gamma} ds, 
\end{align*}
for all $t \in [0,T]$. 
Then, we take the sum of (\ref{conti1}), (\ref{conti2}) \pier{and combine with the above equality. Thanks to}
the monotonicity of $\beta$, $\beta_\Gamma$, 
the {L}ipschitz continuity of $\pi$, $\pi_\Gamma$, 
and the {P}oincar\' e--{W}irtinger inequality, 
we deduce that 
\begin{align*}
	& 
	\bigl| \phi^{(1)}(t)-\phi^{(2)}(t)\bigr|_{V_{0*}}^2
	+
	\bigl| \psi^{(1)}(t)-\psi^{(2)}(t)\bigr|_{V_{\Gamma,0*}}^2
	\nonumber \\
	& \quad {}
	+ 2\int_0^t\!\!\int _\Omega \bigl| \nabla \bigl(\phi^{(1)}-\phi^{(2)} \bigr) \bigr|^2dx ds 
	+ 2\int_0^t\!\!\int _\Gamma \bigl| \nabla_\Gamma \bigl(\psi^{(1)}-\psi^{(2)} \bigr) \bigr|^2d\Gamma ds
	\nonumber \\
	& \le
	\bigl| \phi^{(1)}_0-\phi^{(2)}_0 \bigr|_{V_{0*}}^2
	+
	\bigl| \psi^{(1)}_0-\psi^{(2)}_0 \bigr|_{V_{\Gamma, 0*}}^2 
	+
	2L \int_0^t\!\!\int _\Omega \bigl| \phi^{(1)}-\phi^{(2)}\bigr|^2dx ds 
	\nonumber \\
	& \quad {} 
	+ 
	2L_\Gamma \int_0^t\!\!\int _\Gamma \bigl| \psi^{(1)}-\psi^{(2)} \bigr|^2d\Gamma ds 
	+
	\frac{1}{2} \int_0^t \bigl| \phi^{(1)}-\phi^{(2)}\bigr|^2_{V_0} ds 
	+ \frac{C_{\rm P}}{2} \bigl| f^{(1)}-f^{(2)} \bigr|_{L^2(0,T;V^*)}^2 
	\nonumber \\
	& \quad {}
	+
	\frac{1}{2} \int_0^t \bigl| \psi^{(1)}-\psi^{(2)}\bigr|^2_{V_{\Gamma,0}} ds 
	+ \frac{C_{\rm P}}{2} \bigl| g^{(1)}-g^{(2)} \bigr|_{L^2(0,T;V_{\Gamma}^*)}^2, 
\end{align*}
for all $t \in [0,T]$. 
Here, we observe that 
\begin{align*}
	& 2L \int_0^t\!\!\int _\Omega \bigl| \phi^{(1)}-\phi^{(2)}\bigr|^2dx ds 
	\nonumber \\
	&  = 2L \int_0^t\!\!\int _\Omega \nabla 
	{\mathcal N}_\Omega \bigl( \phi^{(1)}-\phi^{(2)} \bigr) \cdot \nabla \bigl( \phi^{(1)}-\phi^{(2)} \bigr) dx ds 
	\nonumber \\
	& \le 
	\frac{1}{2}\int_0^t\!\!\int _\Omega \bigl| \nabla \bigl(\phi^{(1)}-\phi^{(2)} \bigr) \bigr|^2dx ds 
	+ 2L^2 \int_0^t\!\!\int _\Omega 
	\bigl| \nabla {\mathcal N}_\Omega\bigl(\phi^{(1)}-\phi^{(2)} \bigr) \bigr|^2dx ds 
	\nonumber \\
	&  =
	\frac{1}{2} \int_0^t \bigl| \phi^{(1)}-\phi^{(2)}  \bigr|^2_{V_0} ds 
	+ 2L^2 \int_0^t\bigl|\phi^{(1)}-\phi^{(2)}  \bigr|^2_{V_{0*}} ds, 
\end{align*}
and similarly,
\begin{equation*}
	2L_\Gamma \int_0^t\!\!\int _\Gamma \bigl| \psi^{(1)}-\psi^{(2)}\bigr|^2d\Gamma ds 
	\le \frac{1}{2}
	\int_0^t \bigl| \psi^{(1)}-\psi^{(2)} \bigr|^2_{V_{\Gamma,0}} ds 
	+ 2 L^2_\Gamma \int_0^t\bigl|\psi^{(1)}-\psi^{(2)}  \bigr|^2_{V_{\Gamma,0*}} ds,
\end{equation*}
for all $t \in [0,T]$. Therefore, applying the {G}ronwall lemma and \pier{invoking the} 
equivalences of norms, we conclude the proof of 
Theorem~2.2. As an immediate sequence, the continuous dependence implies the uniqueness of the weak solution obtained in Theorem~2.1. 
\hfill $\Box$\medskip

The continuous dependence \pier{estimate} can be 
extended to the viscous or partially viscous cases, with the following modification: there exists a constant $C$ such that \pier{the inequality}
\begin{align*}
	& 
	\bigl| \phi^{(1)}-\phi^{(2)}
	\bigr|_{C([0,T];V^*)} 
	+
	\bigl| \psi^{(1)}-\psi^{(2)} 
	\bigr|_{C([0,T];V_\Gamma^*)} 
	+
	\bigl| \phi^{(1)}-\phi^{(2)} 
	\bigr|_{L^2(0,T;V)} 
	+
	\bigl| \psi^{(1)}-\psi^{(2)} 
	\bigr|_{L^2(0,T;V_\Gamma)} 
	\nonumber \\
	& \quad {} + 
	\sqrt{\tau} \bigl| \phi^{(1)}-\phi^{(2)}
	\bigr|_{C([0,T];H)} 
	+ \sqrt{\sigma}\bigl| \psi^{(1)}-\psi^{(2)} 
	\bigr|_{C([0,T];H_\Gamma)} 
	\nonumber \\
	& \le C \left\{ 
	\bigl| \phi^{(1)}_0-\phi^{(2)}_0
	\bigr|_{V^*}
	+
	\bigl| \psi^{(1)}_0-\psi^{(2)}_0
	\bigr|_{V_\Gamma^*}
	+
	\bigl| f^{(1)} - f^{(2)}
	\bigr|_{L^2(0,T;V^*)} 
	+ 
	\bigl| g^{(1)}-g^{(2)}
	\bigr|_{L^2(0,T;V_\Gamma^*)}
	\right. 
	\nonumber \\
	& \quad 
	{} + \left.  \sqrt{\tau} \bigl| \phi^{(1)}_0-\phi^{(2)}_0
	\bigr|_{H}
	+
	\sqrt{\sigma}\bigl| \psi^{(1)}_0-\psi^{(2)}_0
	\bigr|_{H_\Gamma}\right\}
\end{align*} 
holds for $\tau \ge 0$ and $\sigma \ge 0$.

\section{Existence of strong solution}
\setcounter{equation}{0}

In this section, we establish a regularity result, which leads to the existence of a strong solution in the case of the pure 
{C}ahn--{H}illiard system \eqref{LW1}--\eqref{LW7}. 

Now, we point out the additional assumptions we need on the given data: 
\begin{enumerate}
 \item[(A6)] $f \in H^1(0,T;H)$ and $g \in H^1(0,T;H_\Gamma)$;
 \item[(A7)] $-\Delta \phi_0 + \beta_\varepsilon (\phi_0) + \pi(\phi_0)-f(0)$ remains bounded in $V$ as $\varepsilon \to 0$, and 
 $\partial _{\boldsymbol{\nu}} \phi_0- \Delta_\Gamma \psi_0 
 + \beta_{\Gamma, \varepsilon} (\psi_0)+\pi_\Gamma(\psi_0)-g(0)$ remains bounded in $V_\Gamma$ as 
 $\varepsilon \to 0$.
\end{enumerate}
\smallskip

Our third result of this paper is related to the existence of strong solutions:\smallskip

\paragraph{\bf Theorem 4.1.} 
{\it Under assumptions {\rm (A1)}--{\rm (A7)}, \pier{the unique weak solution of problem \eqref{LW1}--\eqref{LW7} obtained in Theorem~2.1 is a strong one in the sense of Definition~2.2.}}\smallskip

\paragraph{\bf Proof.} 
Let us take the difference of equations (\ref{td1}) \pier{written for $n$ and} $n-1$. Then multiplying the resultant by 
\begin{equation*}
	z:={\mathcal N}_\Omega( \phi_{n+1}+h \mu_{n+1}-\phi_{n}-h \mu_{n} )
\end{equation*}
and integrating the resultant over $\Omega$, we obtain 
\begin{align*}
	& \frac{1}{h} 
	\int_\Omega \left( \phi_{n+1}+h \mu_{n+1}-2 \phi_{n} - 2 h \mu_{n}+\phi_{n-1} + h \mu_{n-1} \right) z dx 
	\nonumber \\
	& \quad {} + \int_\Omega \nabla (\mu_{n+1}-\mu_n) \cdot \nabla z dx =0,
\end{align*}
whence 
\begin{align}
	& \frac{1}{2h} 
	| \phi_{n+1}+h \mu_{n+1}-\phi_{n} - h \mu_{n} |_{V_{0*}}^2 
	-
	\frac{1}{2h} 
	| \phi_{n}+h \mu_{n}-\phi_{n-1} - h \mu_{n-1} |_{V_{0*}}^2 
	\nonumber \\
	& \quad {} + \frac{1}{2h} 
	| \phi_{n+1}+h \mu_{n+1}-2\phi_n-2h \mu_n + \phi_{n-1} + h \mu_{n-1} |_{V_{0*}}^2 
	\nonumber \\
	& \quad {}+
	\int_\Omega ( \mu_{n+1}-\mu_n ) 
	( \phi_{n+1}+h \mu_{n+1}-\phi_n - h \mu_n ) dx = 0, 
	\label{st1}
\end{align}
for all $n=\pier{1,\ldots,N-1}$.
Next, we use (\ref{td2}) to derive that
\begin{align}
	& 
	\int_\Omega ( \mu_{n+1}-\mu_n ) 
	( \phi_{n+1}+h \mu_{n+1}-\phi_n - h \mu_n) dx 
	\nonumber \\
	&  = \frac{\tau}{2h} | \phi_{n+1} - \phi_n |_H^2 
	- \frac{\tau}{2h} | \phi_n - \phi_{n-1} |_H^2 
	+ \frac{\tau}{2h} | \phi_{n+1} -2 \phi_n + \phi_{n-1} |_H^2 
	\nonumber \\
	& \quad {} 
	+ \int_\Omega \bigl| \nabla (\phi_{n+1}-\phi_n) \bigr|^2 dx 
	- \int_\Gamma \partial _{\boldsymbol{\nu}}  (\phi_{n+1}-\phi_n)  (\psi_{n+1}-\psi_n) d\Gamma
	\nonumber \\
	& \quad {} 
	+ \int_\Omega \bigl( \beta_\varepsilon(\phi_{n+1})-\beta_\varepsilon(\phi_n) \bigr) (\phi_{n+1}-\phi_n) dx 
	+ \int_\Omega \bigl( \pi(\phi_{n+1})-\pi(\phi_n) \bigr) (\phi_{n+1}-\phi_n) dx 
	\nonumber \\
	& \quad  {}
	- \int_\Omega (f_n-f_{n-1}) (\phi_{n+1}-\phi_n) dx 
	+ h \int_\Omega | \mu_{n+1}-\mu_n |^2 dx, 
	\label{st2}
\end{align}
for all $n=0,1,\ldots,N-1$. 
Next, in order to treat the fifth term on the right hand side of (\ref{st2}), we recall equation (\ref{td6}) and infer that 
\begin{align}
	& 
	- \int_\Gamma \partial _{\boldsymbol{\nu}}  (\phi_{n+1}-\phi_n)  (\psi_{n+1}-\psi_n) d\Gamma
	\nonumber \\
	& 
	 = - \int_\Gamma (w_{n+1}-w_n ) 
	( \psi_{n+1}- \psi_n ) d\Gamma 
	+ \frac{\sigma}{2h} | \psi_{n+1}-\psi_n |_{H_\Gamma}^2 
	- \frac{\sigma}{2h} | \psi_n-\psi_{n-1} |_{H_\Gamma}^2 
	\nonumber \\
	& \quad {}
	+ \frac{\sigma}{2h} | \psi_{n+1}-2 \psi_n + \psi_{n-1} |_{H_\Gamma}^2 
	+ \int_\Gamma \bigl| \nabla_\Gamma \bigl( \psi_{n+1}-\psi_n \bigr) \bigr|^2 d\Gamma
	\nonumber \\
	& \quad {}
	+ \int_\Gamma
	\bigl( \beta_{\Gamma,\varepsilon} (\psi_{n+1})-\beta_{\Gamma,\varepsilon}(\psi_n) \bigr) (\psi_{n+1}-\psi_n) d\Gamma \nonumber\\
	&\quad 
	+ \int_\Gamma \bigl( \pi_\Gamma (\psi_{n+1})-\pi_\Gamma (\psi_n) \bigr) (\psi_{n+1}-\psi_n) d\Gamma 
	\nonumber \\
	& \quad {}
	- \int_\Gamma ( g_n-g_{n-1}) (\psi_{n+1}-\psi_n) d\Gamma, 
	\label{st3}
\end{align}
for all $n=0,1,\ldots,N-1$. 
We now exploit (\ref{td5}) to discuss the first term on the right hand side of (\ref{st3}). Taking the test function
\begin{equation*}
	z_\Gamma:={\mathcal N}_\Gamma (\psi_{n+1} + h w_{n+1} -\psi_n - h w_n),
\end{equation*}
\pier{we have that }
\begin{align}
	& \frac{1}{h} \int_\Gamma \bigl( \psi_{n+1} + h w_{n+1} -\psi_n - h w_n 
	- ( \psi_n + h w_n -\psi_{n-1} - h w_{n-1} ) \bigr)  z_\Gamma d\Gamma 
	\nonumber \\
	&\quad  = \bigl( \Delta_\Gamma (w_{n+1}-w_n), z_\Gamma \bigr)_{H_\Gamma} 
	\nonumber \\
	&\quad  = - \int_\Gamma \nabla _\Gamma (w_{n+1}-w_n) \cdot \nabla_\Gamma z_\Gamma d\Gamma 
	\nonumber \\
	&\quad  = - \int_\Gamma 
	 (w_{n+1}-w_n ) ( \psi_{n+1}- \psi_n )d\Gamma 
	- h \int_\Gamma |w_{n+1}-w_n|^2 d\Gamma,
	\label{st4}
\end{align}
for all $n=\pier{1,\ldots,N-1}$.  

Now, by collecting the identities (\ref{st1})--(\ref{st4}) \pier{and using {\rm (A4)}, we deduce from Young's inequality that}
\begin{align*}
	& 
	\frac{1}{2h} | \phi_{n+1}+ h \mu_{n+1}-\phi_n-h \mu_n | _{V_{0*}}^2 
	+ \frac{1}{2h} | \psi_{n+1}+ h w_{n+1}-\psi_n-h w_n | _{V_{\Gamma,0*}}^2 
	\nonumber \\
	& \quad {}
	- \frac{1}{2h} | \phi_{n}+ h \mu_{n}-\phi_{n-1}- h \mu_{n-1} | _{V_{0*}}^2 
	- \frac{1}{2h} | \psi_{n}+ h w_{n}-\psi_{n-1}-h w_{n-1} | _{V_{\Gamma,0*}}^2 
	\nonumber \\
	& \quad {}
	+ \frac{1}{2h} | \phi_{n+1}+ h \mu_{n+1}-2 \phi_n-2 h \mu_n + \phi_{n-1}+h\mu_{n-1}| _{V_{0*}}^2 
	\nonumber \\
	& \quad {} 
	+ \frac{1}{2h} | \psi_{n+1}+ h w_{n+1}-2\psi_n-2h w_n + \psi_{n-1}+h w _{n-1}| _{V_{\Gamma,0*}}^2 
	\nonumber \\
	& \quad {}	
	+ \frac{\tau}{2h} | \phi_{n+1}-\phi_n| _{H}^2 
	+ \frac{\sigma}{2h} | \psi_{n+1}-\psi_n | _{H_\Gamma }^2 
	- \frac{\tau}{2h} | \phi_{n}-\phi_{n-1} | _{H}^2 
	- \frac{\sigma}{2h} | \psi_{n}-\psi_{n-1} | _{H_\Gamma}^2 
	\nonumber \\
	& \quad {}
	+ \frac{\tau}{2h} | \phi_{n+1}-2 \phi_n+ \phi_{n-1}| _{H}^2 
	+ \frac{\sigma}{2h} | \psi_{n+1}-2\psi_n+ \psi_{n-1}| _{H_\Gamma}^2 
	+ \int_\Omega \bigl| \nabla (\phi_{n+1}-\phi_n) \bigr|^2 dx 
	\nonumber \\
	& \quad {}
	+ \int_\Gamma  \bigl| \nabla_\Gamma (\psi_{n+1}-\psi_n) \bigr|^2 d\Gamma 
	+ h \int_\Omega  | \mu_{n+1}-\mu_n |^2 dx 
	+ h \int_\Gamma  | w_{n+1}-w_n |^2 d\Gamma
	\nonumber \\
	& 
	\le (L+1) \int_\Omega 
	| \phi_{n+1}- \phi_n |^2 dx 
	+ \frac{1}{4} \int_\Omega |f_n-f_{n-1}|^2 dx 
	+
	(L_\Gamma+1) 
	\int_\Gamma 
	| \psi_{n+1}- \psi_n |^2 d\Gamma  
	\nonumber \\
	& \quad 
	{} + \frac{1}{4}\int_\Gamma | g_n-g_{n-1}|^2 d\Gamma 
\end{align*}
for $n=1,2,\ldots,N-1$. We divide the above inequality by $h$ and sum up for $n=1$ to $n=m$, by finding that
 \begin{align}
	& 
	\frac{1}{2} \left| \frac{\phi_{m+1}-\phi_m}{h}+ \mu_{m+1}-\mu_m \right| _{V_{0*}}^2 
	+ \frac{1}{2} \left| \frac{\psi_{m+1}-\psi_m}{h}+ w_{m+1}-w_m  \right| _{V_{\Gamma,0*}}^2 
	\nonumber \\
	& \quad {}
	+ \frac{\tau}{2} \left| \frac{\phi_{m+1}-\phi_m}{h}\right| _{H}^2 
	+ \frac{\sigma}{2} \left| \frac{\psi_{m+1}-\psi_m}{h}\right| _{H_\Gamma}^2 
	+ \sum_{n=1}^m h 
	\int_\Omega \left| \nabla \left( \frac{\phi_{n+1}-\phi_n}{h} \right) \right|^2 dx 
	\nonumber \\
	& \quad {}
	+ \sum_{n=1}^m h 
	\int_\Gamma \left| \nabla_\Gamma \left( \frac{\psi_{n+1}-\psi_n}{h} \right) \right|^2 d\Gamma 
	+ \sum_{n=1}^m  | \mu_{n+1}-\mu_n |^2_H
	+ \sum_{n=1}^m  | w_{n+1}-w_n |^2_{H_\Gamma} 
	\nonumber \\
	& 
	\le 
	\frac{1}{2} \left| \frac{\phi_1-\phi_0}{h} + \mu_1-\mu_0 \right| _{V_{0*}}^2 
	+ \frac{1}{2} \left| \frac{\psi_1-\psi_0}{h} + w_1-w_0 \right| _{V_{\Gamma,0*}}^2 
	+ \frac{\tau}{2} \left| \frac{\phi_1-\phi_0}{h}\right| _{H}^2 
	\nonumber \\
	& \quad 
	{} 
	+ \frac{\sigma}{2} \left| \frac{\psi_1-\psi_0}{h}\right| _{H_\Gamma}^2 
	+ (L+1) \sum_{n=1}^m h 
	\left| \frac{\phi_{n+1}- \phi_n}{h} \right|^2_H  
	+\frac{1}{4}
	\sum_{n=1}^m h 
	\left| \frac{f_n- f_{n-1}}{h} \right|^2_H  
	\nonumber \\
	& \quad 
	{} + (L_\Gamma +1 )
	\sum_{n=1}^m h 
	\left| \frac{\psi_{n+1}- \psi_n}{h} \right|^2_{H_\Gamma} 
	+ \frac{1}{4}\sum_{n=1}^m h 
	\left| \frac{g_n- g_{n-1}}{h} \right|^2_{H_\Gamma}, 
	\label{st5}
\end{align}
for all $m=1,2,\ldots,N-1$. In order to estimate the first four terms on the right hand side of \eqref{st5}, 
we multiply (\ref{td1}) at $n=0$ by 
\begin{equation*}
	\frac{1}{h} {\mathcal N}_\Omega ( \phi_1+h \mu_1-\phi_0-h \mu_0 )
\end{equation*}
and obtain 
\begin{equation}
	\left| \frac{\phi_1-\phi_0}{h}+\mu_1-\mu_0 \right|_{V_{0*}}^2 
	+ \int_\Omega \mu_1 \left( \frac{\phi_1-\phi_0}{h}+\mu_1-\mu_0 \right) dx =0.
	\label{st6}
\end{equation}
From (\ref{td2}) \pier{and \eqref{td7}}, it follows that
\begin{align}
	& \int_\Omega \mu_1 \left( \frac{\phi_1-\phi_0}{h}+\mu_1-\mu_0 \right) dx 
	\nonumber \\
	& = \tau \int_\Omega \left| \frac{\phi_1-\phi_0}{h} \right| ^2 dx 
	+ \int_\Omega \nabla \phi_1 \cdot \nabla \left( \frac{\phi_1-\phi_0}{h} \right) dx 
	- \int_\Gamma \partial_{\boldsymbol{\nu}} \phi_1 \left( \frac{\psi_1-\psi_0}{h} \right) d\Gamma
	\nonumber \\
	& \quad {}
	+ \int_\Omega \bigl( \beta_\varepsilon (\phi_1)- \beta_\varepsilon (\phi_0)\bigr) \left( \frac{\phi_1-\phi_0}{h} \right) dx 
	+ \int_\Omega  \beta_\varepsilon (\phi_0) \left( \frac{\phi_1-\phi_0}{h} \right) dx 
	\nonumber \\
	& \quad {}
	+ \int_\Omega \bigl( \pi (\phi_1)- \pi (\phi_0)\bigr) \left( \frac{\phi_1-\phi_0}{h} \right) dx
	+ \int_\Omega  \pi (\phi_0) \left( \frac{\phi_1-\phi_0}{h} \right) dx 
	\nonumber \\
	& \quad {}
	- \int_\Omega f_0 \left( \frac{\phi_1-\phi_0}{h} \right) dx+ \int_\Omega |\mu_1-\mu_0|^2 dx\pier{.}
	\label{st7}
\end{align}
\pier{Next,} we multiply (\ref{td6}) at $n=0$ by 
$(\psi_1-\psi_0)/h$ to obtain
\begin{align}
	& - \int_\Gamma \partial_{\boldsymbol{\nu}} \phi_1 \left( \frac{\psi_1-\psi_0}{h} \right) d\Gamma
	\nonumber \\
	& = 
	- \int_\Gamma w_1 	\left( \frac{\psi_1-\psi_0}{h} \right) d\Gamma
	+ \sigma \int_\Gamma \left| \frac{\psi_1-\psi_0}{h} \right| ^2 d\Gamma 
	+ \int_\Gamma \nabla_\Gamma \psi_1 \cdot \nabla_\Gamma \left( \frac{\psi_1-\psi_0}{h} \right) d\Gamma 
	\nonumber \\
	& \quad {}
	+ \int_\Gamma \bigl( \beta_{\Gamma,\varepsilon} (\psi_1)- \beta_{\Gamma,\varepsilon} (\psi_0)\bigr) 
	\left( \frac{\psi_1-\psi_0}{h} \right) d\Gamma 
	+ \int_\Gamma \beta_{\Gamma,\varepsilon} (\psi_0) \left( \frac{\psi_1-\psi_0}{h} \right) d\Gamma 
	\nonumber \\
	& \quad {}
	+ \int_\Gamma \bigl( \pi_\Gamma (\psi_1)- \pi_\Gamma (\psi_0)\bigr) 
	\left( \frac{\psi_1-\psi_0}{h} \right) d\Gamma
	+ \int_\Gamma \pi_\Gamma (\psi_0) \left( \frac{\psi_1-\psi_0}{h} \right) d\Gamma
	\nonumber \\
	& \quad {}
	- \int_\Gamma g_0
	\left( \frac{\psi_1-\psi_0}{h} \right) d\Gamma.
	\label{st8}
\end{align}
Besides, we go back to (\ref{td5}) at $n=0$, multiply it by
\begin{equation*}
	\frac{1}{h} {\mathcal N}_\Gamma ( \psi_1+h w_1-\psi_0-h w_0 ),
\end{equation*}
and integrate the resultant over $\Gamma$, to obtain
\begin{equation}
	\left| \frac{\psi_1-\psi_0}{h}+w_1-w_0 \right|_{V_{\Gamma,0*}}^2 
	+ \int_\Gamma w_1 \left( \frac{\psi_1-\psi_0}{h}+w_1-w_0 \right) d\Gamma =0. 
	\label{st9}
\end{equation}
Collecting (\ref{st6})--(\ref{st9}) and applying the assumption (A7) \pier{along}
with the monotonicity 
of $\beta_\varepsilon$, $\beta_{\Gamma,\varepsilon}$,  
the {L}ipschitz continuity of $\pi$, $\pi_\Gamma$ and Young's inequality,  
we deduce that there exists a positive constant 
$\tilde{M}_9$, independent of $h\in(0,h^{**}]$, $\tau, \sigma, \varepsilon \in (0,1]$, \pier{%
%indeed it is surely independent of v$\varepsilon \in (0,1]$ from (A7), 
such that}
\begin{align}
	& \left| \frac{\phi_1-\phi_0}{h}+\mu_1-\mu_0 \right|_{V_{0*}}^2 +
	\left| \frac{\psi_1-\psi_0}{h}+w_1-w_0 \right|_{V_{\Gamma,0*}}^2 
	+\tau \int_\Omega \left| \frac{\phi_1-\phi_0}{h} \right| ^2 dx 
	\nonumber \\
	& \quad {} +\sigma \int_\Gamma \left| \frac{\psi_1-\psi_0}{h} \right| ^2 d\Gamma 
	+ h \int_\Omega \left| \nabla \left( \frac{\phi_1-\phi_0}{h} \right) \right|^2 dx 
	+ h \int_\Gamma \left| \nabla_\Gamma \left( \frac{\psi_1-\psi_0}{h} \right) \right|^2 d\Gamma 
	\nonumber \\
	& \quad {} + \int_\Omega |\mu_1-\mu_0|^2 dx + \int_\Gamma |w_1-w_0|^2 d\Gamma
	\nonumber \\
	& \le -
	\int_\Omega \bigl( -\Delta \phi_0 + \beta_\varepsilon (\phi_0)+ \pi (\phi_0) - f_0 \bigr) 
	\left( \frac{\phi_1-\phi_0}{h} \right) dx 
	+ L  h\left| \frac{\phi_1-\phi_0}{h} \right|_H^2 
	\nonumber \\
	& \quad {}
	- \int_\Gamma \bigl( \partial_{\boldsymbol{\nu}} \phi_0 
	- \Delta_\Gamma \psi_0 + \beta_{\Gamma,\varepsilon}(\psi_0)+\pi_\Gamma(\psi_0)-g_0
	\bigr)
	\left( \frac{\psi_1-\psi_0}{h} \right) d\Gamma
	+ L_\Gamma h \left| \frac{\psi_1-\psi_0}{h} \right|_{H_\Gamma}^2
	\nonumber \\
	& \le 
	\tilde{M}_9 
	+ \frac{1}{2}\left| \frac{\phi_1-\phi_0}{h}+\mu_1-\mu_0 \right|_{V_{0*}}^2
	+ \frac{1}{2} |\mu_1-\mu_0|^2 _H
	+ L h\left| \frac{\phi_1-\phi_0}{h} \right|_H^2 
	\nonumber \\
	& \quad {}
	+ \frac{1}{2} \left| \frac{\psi_1-\psi_0}{h} +w_1-w_0\right|_{V_{\Gamma,0*}}^2 
	+ \frac{1}{2} |w_1-w_0|^2 _{H_\Gamma}
	+ L_\Gamma h \left| \frac{\psi_1-\psi_0}{h} \right|_{H_\Gamma}^2.
	\label{st10}
\end{align}
Then, we can add (\ref{st5}) and (\ref{st10}) to obtain that 
\begin{align*}
	& 
	\frac{1}{2} \left| \frac{\phi_{m+1}-\phi_m}{h}+ \mu_{m+1}-\mu_m \right| _{V_{0*}}^2 
	+ \frac{1}{2} \left| \frac{\psi_{m+1}-\psi_m}{h}+ w_{m+1}-w_m  \right| _{V_{\Gamma,0*}}^2 
	\nonumber \\
	& \quad {}
	+ \frac{\tau}{2} \left| \frac{\phi_{m+1}-\phi_m}{h}\right| _{H}^2 
	+ \frac{\sigma}{2} \left| \frac{\psi_{m+1}-\psi_m}{h}\right| _{H_\Gamma}^2 
	+ \sum_{n=0}^m h 
	\int_\Omega \left| \nabla \left( \frac{\phi_{n+1}-\phi_n}{h} \right) \right|^2 dx 
	\nonumber \\
	& \quad {}
	+ \sum_{n=0}^m h 
	\int_\Gamma \left| \nabla_\Gamma \left( \frac{\psi_{n+1}-\psi_n}{h} \right) \right|^2 d\Gamma 
	+ \frac{1}{2} \sum_{n=0}^m  | \mu_{n+1}-\mu_n |^2_H
	+ \frac{1}{2}\sum_{n=0}^m  | w_{n+1}-w_n |^2_{H_\Gamma} 
	\nonumber \\
	& 
	\le 
	M_9 
	+ (L +1) \sum_{n=0}^m h\left| \frac{\phi_{n+1}-\phi_n}{h} \right|_H^2 
	+ (L_\Gamma+1) \sum_{n=0}^m h \left| \frac{\psi_{n+1}-\psi_n}{h} \right|_{H_\Gamma}^2.
\end{align*}
for all $m=1,2,\ldots,N-1$, where $M_9$ is a positive constant \pier{depending} on 
$|\partial _t f|_{L^2(0,T;H)}$, $|\partial_t g|_{L^2(0,T;H_\Gamma)}$ and $\tilde{M}_9$.  
Finally, we apply the compactness inequality 
on the right hand side of that above inequality as follows:
\begin{align*}
	& (L +1) \sum_{n=0}^m h\left| \frac{\phi_{n+1}-\phi_n}{h} \right|_H^2 
	+
	(L_\Gamma+1) \sum_{n=0}^m h \left| \frac{\psi_{n+1}-\psi_n}{h} \right|_{H_\Gamma}^2
	\nonumber \\
	& \le \delta \sum_{n=0}^m h \int_\Omega \left| \nabla \left( \frac{\phi_{n+1}-\phi_n}{h} \right) \right|^2 dx 
	+ \delta \sum_{n=0}^m h 
	\int_\Gamma \left| \nabla_\Gamma \left( \frac{\psi_{n+1}-\psi_n}{h} \right) \right|^2 d\Gamma 
	\nonumber \\
	& \quad {}
	+ C_\delta \sum_{n=0}^m h \left| \frac{\phi_{n+1}-\phi_n}{h} \right|^2_{V^*}
	+ C_\delta \sum_{n=0}^m h \left| \frac{\psi_{n+1}-\psi_n}{h} \right|^2_{V_\Gamma^*},
\end{align*}
for all $\delta>0$. Observe now that by taking $\delta:=1/2$, 
the last two terms are already bounded due to (\ref{M1}). 
Thus, we conclude that 
there exists a positive constant $M_{10}$, 
independent of $h\in(0,h^{**}]$, $\tau, \sigma, \varepsilon \in (0,1]$, 
such that
\begin{align}
	& \bigl|\partial_t \hat{\phi}_h + h \partial _t \hat{\mu}_h \bigr|_{L^\infty (0,T;V^*)}^2 
	+ \bigl|\partial_t \hat{\psi}_h + h \partial _t \hat{w}_h \bigr|_{L^\infty(0,T;{V_\Gamma}^*)}^2 + \tau \bigl|\partial_t \hat{\phi}_h \bigr|_{L^\infty(0,T;H)}^2 
	\nonumber \\
	& \quad 
	+ \sigma \bigl|\partial_t \hat{\psi}_h \bigr|_{L^\infty (0,T;H_\Gamma)}^2 
	+ \bigl| \partial_t \hat{\phi}_h \bigr|_{L^2(0,T;V)}^2 
	+ \bigl| \partial_t \hat{\psi}_h \bigr|_{L^2(0,T;V_\Gamma)}^2 
	+ h \bigl| \partial_t \hat{\mu}_h\bigr|_{L^2(0,T;H)}^2
	\nonumber \\
	& \quad 
	+ h \bigl| \partial_t \hat{w}_h\bigr|_{L^2(0,T;H_\Gamma)}^2 
	+ \bigl| h \partial_t \hat{\mu}_h\bigr|_{L^\infty (0,T;H)}^2 
	+ \bigl| h \partial_t \hat{w}_h\bigr|_{L^\infty (0,T;H_\Gamma)}^2 
	\le M_{10},
	\label{M10}
\end{align}
for all $h\in(0,h^{**}]$. 
The subsequent estimates repeat the previous \pier{ones}, that is, from Lemmas~3.2 to 3.6 as follows:
\begin{list}{$\vartriangleright$}{}
\item From (\ref{L1}), using (\ref{M10}) we infer that
\begin{equation*}
	\bigl| \bar{\mu}_h-m_\Omega ( \bar{\mu}_h ) \bigr|_{L^\infty(0,T;V)} 
	+ 
	\bigl| \bar{w}_h-m_\Gamma ( \bar{w}_h ) \bigr|_{L^\infty(0,T;V_\Gamma)} 
	 \le \Lambda_{1}(t),
\end{equation*}
for all $h \in(0,h^{**}]$ and $\varepsilon \in (0,1]$, 
with $\Lambda_1$ being defined by (\ref{Lam1}). Now $\Lambda_{1}$ is bounded in $L^\infty(0,T)$;
\item arguing as in (\ref{Lem32-a})--(\ref{Lem32-c}) and checking the right hand side of (\ref{Lem32-c}), 
using (\ref{M1}) and (\ref{M10})
we arrive at (\ref{L2}) with $\Lambda_2 \in L^\infty(0,T)$;
\item we can repeat the estimates in Lemmas~3.3 and 3.4, \pier{so that we arrive 
at (\ref{L3})--(\ref{L5})} with 
$\Lambda_3$, $\Lambda_4$, $\Lambda_5 \in L^\infty(0,T)$; 
\item now, we consider the same elliptic system (\ref{el}) and observe that 
we can derive (\ref{L6}) with $\Lambda_6 \in L^\infty(0,T)$; 
\item instead of (\ref{M6}) and (\ref{M7}), \pier{here we derive} from the above modifications the final estimates
\begin{align*}
	& | \bar{\mu}_h |_{L^\infty (0,T;V)} 
	+ 
	| \bar{w}_h |_{L^\infty (0,T;V_\Gamma)} 
	+
	\bigl| \beta_\varepsilon \bigl( \bar \phi_h \bigr) \bigr|_{L^\infty (0,T;H)}
	+
	\bigl| \beta_{\Gamma, \varepsilon} \bigl( \bar \psi_h \bigr) \bigr|_{L^\infty (0,T;H_\Gamma)}
	\nonumber \\
	&\quad 
	 + \bigl| \bar \phi_h \bigr|_{L^\infty (0,T;H^2(\Omega))}
	+
	\bigl| \bar \psi_h \bigr|_{L^\infty (0,T;H^2(\Gamma))}
	+
	\bigl| \partial _{\boldsymbol{\nu}} \bar \phi _h \bigr|_{L^\infty (0,T;H_\Gamma)}
	\le M_{11}, 
	\nonumber \\
	& \sqrt{\tau} | \bar \mu _h |_{L^\infty (0,T;W)} 
	+ \sqrt{\sigma} | \bar w _h |_{L^\infty (0,T;H^2(\Gamma))} \le M_{11},
\end{align*}
for all $h \in(0,h^{**}]$ and $\varepsilon \in (0,1]$, 
where $M_{11}$ is a positive constant
independent of $h\in(0,h^{**}]$, $\tau, \sigma, \varepsilon \in (0,1]$.
\end{list}
Thus, we have obtained sufficient additional estimates that can be extended to the limit functions 
as $h_k \to 0$, by weak or weak star lower semicontinuity of norms. 
In particular, the approximate solution $(\phi_\varepsilon, \mu_\varepsilon, \psi_\varepsilon, w_\varepsilon)$ to problem 
(\ref{yLW1})--(\ref{yLW7}), 
which is unique due to the continuous dependence estimate \pier{stated in} Theorem~2.2, 
additionally satisfies (cf. (\ref{est1})--(\ref{est15}))
\begin{align*}
	&| \phi_\varepsilon |_{W^{1,\infty} (0,T;V^*)} 
	+ 
	| \phi_\varepsilon |_{H^1(0,T;V)} 
	+
	| \phi_\varepsilon |_{L^\infty(0,T;H^2(\Omega))} 
	+
	\sqrt{\tau}
	| \phi_\varepsilon |_{W^{1,\infty}(0,T;H)} 
	\nonumber \\
	&\quad 
	+
	| \mu _\varepsilon |_{L^\infty (0,T;V)}
	+
	| \mu _\varepsilon |_{L^2 (0,T;W \cap H^3(\Omega))}
	+
	\sqrt{\tau} 
	| \mu _\varepsilon |_{L^\infty (0,T;W)}
	\nonumber \\
	&\quad 
	+ | \psi_\varepsilon |_{W^{1,\infty} (0,T;V_\Gamma^*)} 
	+ 
	| \psi_\varepsilon |_{H^1(0,T;V_\Gamma)} 
	+
	| \psi_\varepsilon |_{L^\infty(0,T;H^2(\Gamma))} 
	+
	\sqrt{\sigma}
	| \psi_\varepsilon |_{W^{1,\infty}(0,T;H_\Gamma)} 
	\nonumber \\
	&\quad 
	+
	| w_\varepsilon |_{L^\infty (0,T;V_\Gamma)}
	+
	| w_\varepsilon |_{L^2 (0,T;H^3(\Gamma))}
	+
	\sqrt{\sigma} 
	| w_\varepsilon |_{L^\infty (0,T;H^2(\Gamma))}
	\nonumber \\
	&\quad 
	+
	| \partial _{\boldsymbol{\nu}} \phi _\varepsilon |_{L^\infty (0,T;H_\Gamma)}
	+\bigl| \beta_\varepsilon (\phi_\varepsilon) \bigr|_{L^\infty (0,T;H)}
	+
	\bigl| \beta_{\Gamma, \varepsilon}(\psi_\varepsilon ) \bigr|_{L^\infty (0,T;H_\Gamma)}
	\le M_{12},
\end{align*}
for all $\varepsilon \in (0,1]$, 
where $M_{12}$ is a positive constant
independent of $\tau, \sigma, \varepsilon \in (0,1]$. 
Let us point out that, in order to obtain the above estimate, we have to 
make comparison of terms in (\ref{yLW1}) with (\ref{yLW3}) and (\ref{yLW5}), 
as well as to apply the elliptic regularity results.  
In this framework, the final regularity of the sextuplet  
$(\phi, \mu, \xi, \psi, w, \zeta)$ solving problem (\ref{LW1})--(\ref{LW7}) is
\begin{align*}
	& \phi \in W^{1,\infty} (0,T;V^*)\cap H^1 (0,T;V) \cap L^\infty \bigl( 0,T;H^2(\Omega) \bigr), \\
	& \mu \in L^\infty (0,T;V) \cap L^2 \bigl(0,T;W\cap H^3(\Omega) \bigr), \quad \xi \in L^\infty (0,T;H), \\
	& \psi \in W^{1,\infty} (0,T;V_\Gamma^*)\cap H^1 (0,T;V_\Gamma) \cap L^\infty \bigl(0,T;H^2(\Gamma) \bigr), \\
	& w \in L^\infty (0,T;V_\Gamma) \cap L^2 \bigl(0,T;H^3(\Gamma) \bigr), \quad \zeta \in L^\infty (0,T;H_\Gamma)
\end{align*}
and under these regularities, (\ref{weak1}) and (\ref{weak3}) can be replaced by 
\begin{align*} 
	& \partial_t \phi - \Delta \mu =0 
	& \mbox{a.e.\ in }Q, 
	\\
	& \partial_{\boldsymbol{\nu}} \mu =0 
	& \mbox{a.e.\ on }\Sigma, 
	\\
	& \partial_t \psi -\Delta_\Gamma w = 0 
	& \mbox{a.e.\ on }\Sigma.
\end{align*}
Moreover, in order to obtain 
the regularities for $\mu$ and $w$ stated above,   
we \pier{simply use the weak star lower semicontinuity property}.
Thus, we \pier{complete the proof of Theorem~4.1}. 
\hfill $\Box$\medskip

%%%%%%%%%%%%%%%%%%%%%%%%%%%%%%

\section*{Acknowledgments}
\pcol{This work was started during a visit of P.~Colli and T.~Fukao at the 
School of Mathematical Sciences of Fudan University: the contributed support and 
warm hospitality of the School are very gratefully acknowledged.
P.~Colli also acknowledges other} \pier{support from the Italian Ministry of Education, 
University and Research~(MIUR): Dipartimenti di Eccellenza Program (2018--2022) 
-- Dept.~of Mathematics ``F.~Casorati'', University of Pavia, \pcol{and the} 
GNAMPA (Gruppo Nazionale per l'Analisi Matematica, 
la Probabilit\`a e le loro Applicazioni) of INdAM (Isti\-tuto 
Nazionale di Alta Matematica).}
\takeshi{T.~Fukao acknowledges the support
from the JSPS KAKENHI Grant-in-Aid for Scientific Research(C), Japan Grant Number 17K05321.}
\hao{H.~Wu acknowledges the support from the NNSFC Grant No. 11631011 and the Shanghai Center for Mathematical Sciences at Fudan University.}

\end{document}